# ON 2D NEWEST VERTEX BISECTION: OPTIMALITY OF MESH-CLOSURE AND $H^1$-STABILITY OF $L_2$-PROJECTION

MICHAEL KARKULIK, DAVID PAVLICEK, AND DIRK PRAETORIUS

*Dedicated to Carsten Carstensen on the occasion of his 50th birthday.*

Abstract. Newest vertex bisection (NVB) is a popular local mesh-refinement strategy for regular triangulations which consist of simplices. For the 2D case, we prove that the mesh-closure step of NVB, which preserves regularity of the triangulation, is quasi-optimal and that the corresponding $L^2$-projection onto lowest-order Courant finite elements (P1-FEM) is always $H^1$-stable. Throughout, no additional assumptions on the initial triangulation are imposed. Our analysis thus improves results of Binev, Dahmen & DeVore (Numer. Math. **97**, 2004), Carstensen (Constr. Approx. **20**, 2004), and Stevenson (Math. Comp. **77**, 2008) in the sense that all assumptions of their theorems are removed. Consequently, our results relax the requirements under which adaptive finite element schemes can be mathematically guaranteed to convergence with quasi-optimal rates.

**Corresponding Author:**

Michael Karkulik
Vienna University of Technology
Institute for Analysis and Scientific Computing
Wiedner Hauptstraße 8–10
1040 Wien, Austria

michael.karkulik@tuwien.ac.at
http://www.asc.tuwien.ac.at/∼mkarkulik/

**Acknowledgement:**

The research of the authors M. Karkulik and D. Praetorius is supported through the FWF project *Adaptive Boundary Element Method*, funded by the Austrian Science Fund (FWF) under grant P21732.

http://www.asc.tuwien.ac.at/abem

## 1. Introduction

Contemporary proofs of convergence rates of adaptive finite element methods, e.g., [3, 8, 23], and adaptive boundary element methods [11, 12] rely, among others, on certain properties of the mesh refinement used. The current analysis is mostly based on the newest vertex bisection algorithm (NVB). To be more precise, two properties of the mesh-refinement are currently used: First, one has to relate the optimal mesh $\mathcal{T}_\star$ which could be obtained from the initial mesh $\mathcal{T}_0$ by NVB, with the deterministic mesh $\mathcal{T}_\ell$ which is generated by the adaptive algorithm. Since no relation of these meshes is known a priori, one considers the coarsest common refinement $\mathcal{T}_\star \oplus \mathcal{T}_\ell$ of both meshes. Then, the analysis requires that $\mathcal{T}_\star \oplus \mathcal{T}_\ell$ is not too large. For NVB, it has first been proved in [23] for 2D and later in [8] for $\mathbb{R}^d$, $d \geq 2$, that $\mathcal{T}_\star \oplus \mathcal{T}_\ell$ is, in fact, the overlay of both meshes and that

$$\#\mathcal{T}_\star \oplus \mathcal{T}_\ell \leq \#\mathcal{T}_\star + \#\mathcal{T}_\ell - \#\mathcal{T}_0, \tag{1}$$

independently of $\mathcal{T}_0$. Second, the quasi-optimality analysis relies on the so-called closure estimate of the mesh-refinement: Based on the discrete solution for the mesh $\mathcal{T}_\ell$, an adaptive algorithm chooses a set of marked elements $\mathcal{M}_\ell \subseteq \mathcal{T}_\ell$ which are refined. Refining only marked elements usually leads to hanging nodes, and so the set of marked elements is enlarged to a set of refined elements $\mathcal{R}_\ell \supseteq \mathcal{M}_\ell$, whose refinement yields a regular mesh $\mathcal{T}_{\ell+1}$. Note that refinement of an element $T$ means that $T$ is split into at least two sons, i.e. $\#\mathcal{T}_{\ell+1} - \#\mathcal{T}_\ell \geq \#\mathcal{M}_\ell$ for all $\ell \in \mathbb{N}$. By induction, this provides the bound

$$\sum_{j=0}^{\ell-1} \#\mathcal{M}_j \leq \#\mathcal{T}_\ell - \#\mathcal{T}_0 \quad \text{for all } \ell \in \mathbb{N}.$$

Now, the closure estimate states the converse estimate up to some multiplicative constant

$$\#\mathcal{T}_\ell - \#\mathcal{T}_0 \leq C_1 \sum_{j=0}^{\ell-1} \#\mathcal{M}_j \quad \text{for all } \ell \in \mathbb{N}. \tag{2}$$

Moreover, we refer to [15] for an example that an estimate of the type $\#\mathcal{T}_\ell - \#\mathcal{T}_{\ell-1} \leq C\,\#\mathcal{M}_{\ell-1}$ with an $\ell$-independent constant $C > 0$ cannot hold in general.

The closure estimate (2) has first been proved in [3] for NVB in 2D. The proof has later been generalized to $\mathbb{R}^d$ and $d \geq 2$ in [24]. So far, the proofs of [3, 24] rely on an additional BDD-assumption on the initial $\mathcal{T}_0$, which is an appropriate labeling of the edges of $\mathcal{T}_0$ introduced in detail below. However, no effective algorithm is known which guarantees the BDD-assumption.

To satisfy the BDD-assumption in practice, one can consider the uniform bisec(3)-refinement $\widetilde{\mathcal{T}}_0$ of $\mathcal{T}_0$. Then, the latter can easily be guaranteed for $\widetilde{\mathcal{T}}_0$. This initialization step, however, increases the number of elements by a factor 4 in 2D and even a factor 16 in 3D. This is unacceptable for (dense) boundary element matrices and also unattractive for (sparse) finite element matrices, since adaptivity proves to be quasi-optimal in practice — independently of the initial mesh $\mathcal{T}_0$. Our work now gives, at least for 2D FEM resp. 3D BEM, a mathematical answer for this so far only empirical observation in the sense that we completely avoid any assumption on the initial mesh $\mathcal{T}_0$ to prove (2).

Another important question which is treated in this work, concerns the $L_2$-orthogonal projection $\Pi_\ell$ onto the space of piecewise affine, globally continuous functions $\mathcal{S}^1(\mathcal{T}_\ell)$ on a mesh $\mathcal{T}_\ell$. Since $\mathcal{S}^1(\mathcal{T}_\ell) \subset H^s(\Omega)$ for all $0 \le s \le 1$, norm equivalence on finite-dimensional spaces ensures

$$\|\Pi_\ell u\|_{H^s(\Omega)} \le C_2 \|u\|_{H^s(\Omega)}, \tag{3}$$

where $C_2 > 0$ a-priori depends on $\mathcal{T}_\ell$ and $C_2 = 1$ for $s = 0$. Independence of $C_2$ of $\mathcal{T}_\ell$, i.e. $H^s$-stability of $\Pi_\ell$, is of general interest. Possible applications include hybrid coupled finite element domain decomposition methods [1], hybrid boundary element methods [20], and the construction of efficient preconditioners for finite and boundary element methods [22], a-posteriori error estimation in boundary element methods [18], as well as proving convergence of data-perturbed adaptive boundary element methods [13] and quasi-optimal rates for adaptive FEM with inhomogeneous Dirichlet data [2].

The $H^1$-stability of $\Pi_\ell$ is the subject of a number of works: In [10], stability (3) is shown under certain global conditions on the underlying meshes. The work [5] considers a sequence of (globally) quasi-uniform meshes and shows (3) for $s = 1$. In [4], stability (3) for $s = 1$ is shown for meshes with a controlled local change of the elements' area. In [21], the result of [4] is extended to $0 \le s \le 1$. The conditions in [10] and [4] have finally been merged in [6]. In [7], stability (3) for $0 \le s \le 1$ was shown for a sequence of adaptively generated meshes in 2D. Different to the prior works, no growth conditions for the local mesh-size or element areas were imposed. However, again the initial mesh $\mathcal{T}_0$ has to satisfy a certain (weakened) BDD-type assumption.

The second contribution of this work is that we prove $H^1$-stability of $\Pi_\ell$ on adaptive meshes generated by NVB without any conditions on the initial mesh $\mathcal{T}_0$. In particular, (3) holds for all $0 \le s \le 1$, and $C_2 > 0$ depends only on $\mathcal{T}_0$.

Outline. The remainder of this work is organized as follows: Section 2 formulates the NVB algorithm in 2D as well as a generalization which also covers Carstensen's modification of [7], provides precise statements of our main results, and includes a more detailed discussion on the improvements of the state of the art. In Section 3, we collect general notations as well as auxiliary results. The proof of the closure estimate is given in Section 4. The final Section 5 is concerned with the $H^1$-stability of the $L^2$-projection.

Throughout the proofs, the symbols $\lesssim$ and $\gtrsim$ abbreviate $\le$ and $\ge$ up to some generic constant $C > 0$. The symbol $\simeq$ abbreviates that both estimates $\lesssim$ and $\gtrsim$ hold. In all statements, however, the constants as well as their dependence is explicitely given.

## 2. Newest Vertex Bisection Algorithm & Main Results

Let $\Omega$ be either a bounded domain in $\mathbb{R}^2$ or a part of the boundary of a bounded domain in $\mathbb{R}^3$. We say that $\mathcal{T}$ is a ***mesh*** on $\Omega$ if the following three properties hold:

- $\mathcal{T}$ is a finite set of flat, compact, and non-degenerate triangles $T \in \mathcal{T}$,
- $\overline{\Omega}$ is covered by $\mathcal{T}$, i.e. $\overline{\Omega} = \bigcup_{T\in\mathcal{T}} T$,
- for $T_1, T_2 \in \mathcal{T}$ with $T_1 \neq T_2$, the intersection $T_1 \cap T_2$ is either empty or a node or an edge of both, $T_1$ and $T_2$.

Put differently, $\mathcal{T}$ is a regular triangulation in the sense of Ciarlet, and the last point excludes hanging nodes. Clearly, the shape of $\Omega$ is implicitly assumed to be polygonal and piecewise flat since otherwise the second property cannot be satisfied.

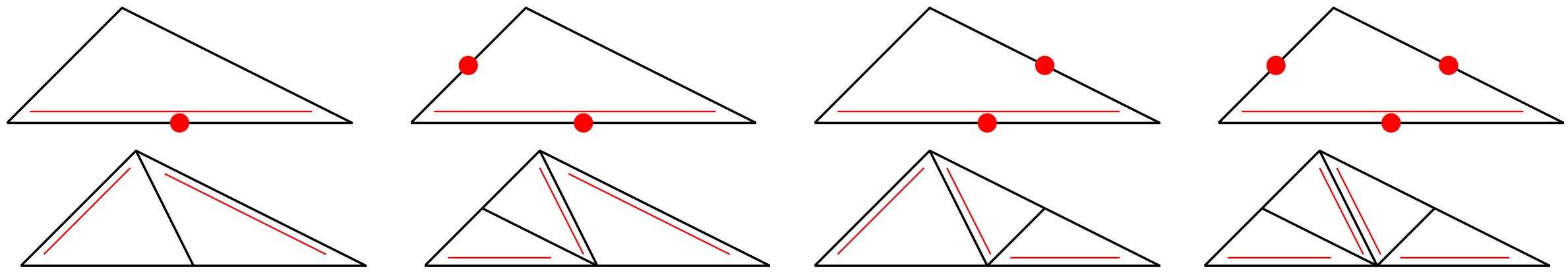

FIGURE 1. For each triangle $T \in \mathcal{T}_\ell$, there is one fixed *reference edge*, indicated by the double line (left, top). Refinement of $T$ is done by bisecting the reference edge, where its midpoint becomes a new node. The reference edges of the son triangles $T' \in \mathcal{T}_{\ell+1}$ are opposite to this newest vertex (left, bottom). To avoid hanging nodes, one proceeds as follows: We assume that certain edges of $T$, but at least the reference edge, are marked for refinement (top). Using iterated newest vertex bisection, the element is then split into 2, 3, or 4 son triangles (bottom). If all elements are refined by three bisections (right, bottom), we obtain the so-called uniform bisec(3)-refinement which is denoted by $\widehat{\mathcal{T}}_\ell$.

**2.1. Newest vertex bisection algorithm.** We briefly describe the idea of NVB: To that end, let $\mathcal{T}_0$ be a given initial mesh. For each triangle $T \in \mathcal{T}_0$, one chooses a so-called *reference edge*, e.g., the longest edge. For NVB, the (inductive) refinement rule reads as follows, where $\mathcal{T}_\ell$ is a regular triangulation already obtained from $\mathcal{T}_0$ by some successive newest vertex bisections:

- To refine an element $T \in \mathcal{T}_\ell$, the midpoint $x_T$ of the reference edge $E_T$ becomes a new node, and $T$ is bisected along $x_T$ and the node opposite to $E_T$ into two sons $T_1$ and $T_2$, see Figure 1 (left).
- As is also shown in Figure 1, the edges opposite to the *newest vertex* $x_T$ become the reference edges of the two son triangles $T_1$ and $T_2$.
- Having bisected all marked triangles, the resulting partition usually has hanging nodes. Therefore, certain additional bisections finally lead to a mesh $\mathcal{T}_{\ell+1}$.

A moment's reflection shows that the latter so-called *mesh-closure* step, which provides the mesh $\mathcal{T}_{\ell+1}$, only leads to finitely many additional bisections. An easy explanation might be the following, which is also illustrated in Figure 1:

- Instead of marked elements, one might think of marked edges.
- If any edge of a triangle $T$ is marked for refinement, we ensure that its reference edge is also marked for refinement. This is done recursively in at most $3 \times \#\mathcal{T}_\ell$ recursions since then all edges are marked for refinement.
- If an element $T$ is bisected, only the reference edge is halved, whereas the other two edges become the reference edges of the two son triangles. The refinement of $T$ into 2, 3, or 4 sons can then be done in one step.

Formally, the algorithm reads as follows:

**Algorithm 1.** INPUT: *given mesh* $\big(\mathcal{T}_\ell, (E_T)_{T\in\mathcal{T}_\ell}\big)$ *and set of marked elements* $\mathcal{M}_\ell \subseteq \mathcal{T}_\ell$.
OUTPUT: *refined mesh* $\big(\mathcal{T}_{\ell+1}, (E_T)_{T\in\mathcal{T}_{\ell+1}}\big)$.

(o) *Set counter* $k := 0$ *and define set of marked reference edges* $\mathcal{M}_\ell^{(0)} := \big\{E_T \,:\, T \in \mathcal{M}_\ell\big\}$.
(i) *Define* $\mathcal{M}_\ell^{(k+1)} := \big\{E_T \,:\, T \in \mathcal{T}_\ell$ *s.t. exists* $E$ *with* $E \in \mathcal{M}_\ell^{(k)}$ *and* $E \subset T\big\}$.

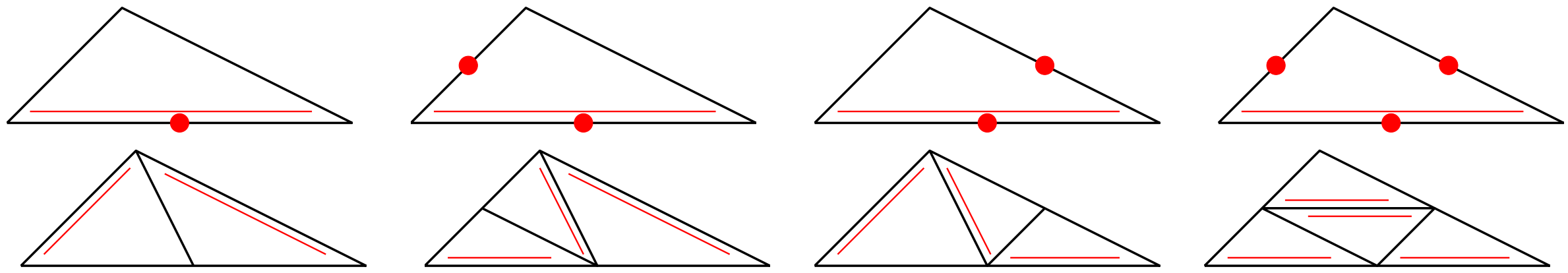

FIGURE 2. The only difference between the refinement rules of the usual newest vertex bisection algorithm (NVB) and CARSTENSEN's modification from [7] is that triangles with three marked edges (right, top) are refined into four similar sons (right, bottom), which is called *red refinement*.

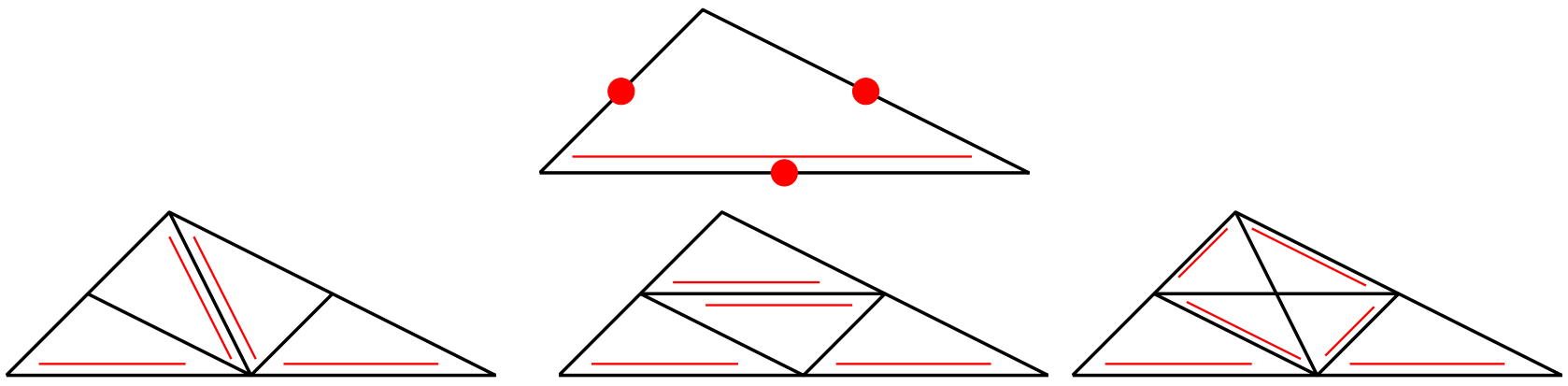

FIGURE 3. The modified newest vertex bisection of Algorithm 2 allows different refinement of an element if all of its edges are marked (top), namely *bisec(3)-refinement* (left), *red refinement* (middle), and *bisec(5)-refinement* (right).

(ii) *If $\mathcal{M}_\ell^{(k)} \subsetneqq \mathcal{M}_\ell^{(k+1)}$, increase counter $k \mapsto k+1$ and goto (i).*

(iii) *Otherwise and with $\mathcal{M}_\ell^{(k)}$ being the set of marked edges, refine each element $\mathcal{T} \in \mathcal{T}_\ell$ according to the rules depicted in Figure 1.*

From now, we identify a mesh $\mathcal{T}$ with the pair $\big(\mathcal{T}, (E_T)_{T\in\mathcal{T}}\big)$, where $E_T$ denotes the reference edge of an element $T \in \mathcal{T}$. If $\mathcal{T}_\ell$ is a mesh and $\mathcal{M}_\ell \subset \mathcal{T}_\ell$ is a set of marked elements, we write $\mathcal{T}_{\ell+1} = \texttt{refineNVB}(\mathcal{T}_\ell, \mathcal{M}_\ell)$ if $\mathcal{T}_{\ell+1}$ is generated by Algorithm 1.

**2.2. Modified newest vertex bisection algorithm.** We now consider a modified (and also generalized) newest vertex bisection (MNVB) algorithm: For each marked element $T \in \mathcal{M}_\ell$, we prescribe (formally even separately) which edges will at least be refined. Then, we proceed analogously to the normal NVB algorithm to determine the set of edges which will be split to avoid hanging nodes. For elements for which one or two edges are marked, we employ the NVB refinement rules of Figure 1. For elements for which all three edges are marked, we assume that additional information is provided to decide whether the element is refined by three bisections, by red-refinement, or by five bisections, see Figure 3. The formal MNVB algorithm reads as follows:

**Algorithm 2.** INPUT: *given mesh $\big(\mathcal{T}_\ell, (E_T)_{T\in\mathcal{T}_\ell}\big)$ and set of marked elements $\mathcal{M}_\ell \subseteq \mathcal{T}_\ell$.*
OUTPUT: *refined mesh $\big(\mathcal{T}_{\ell+1}, (E_T)_{T\in\mathcal{T}_{\ell+1}}\big)$, set of refined elements $\mathcal{R}_\ell := \mathcal{T}_\ell \setminus \mathcal{T}_{\ell+1}$*

(o) *Set counter $k := 0$ and generate set of marked edges $\mathcal{M}_\ell^{(0)} := \big\{E \in \mathcal{E}_\ell \,:\, \text{Exists } T \in \mathcal{M}_\ell \text{ s.t. } E \subset T \text{ has to be refined}\big\}$.*

(i) *Define $\mathcal{M}_\ell^{(k+1)} := \mathcal{M}_\ell^{(k)} \cup \big\{E_T \,:\, T \in \mathcal{T}_\ell \text{ s.t. exists } E \in \mathcal{M}_\ell^{(k)} \text{ with } E \subset T\big\}$.*

(ii) *If* $\mathcal{M}_\ell^{(k)} \subsetneqq \mathcal{M}_\ell^{(k+1)}$*, increase counter* $k \mapsto k+1$ *and goto (i).*

(iii) *Otherwise and with* $\mathcal{M}_\ell^{(k)}$ *being the marked edges, refine each element* $T \in \mathcal{T}_\ell$ *according to the rules depicted in Figure 1 if up to two edges are marked and according to one of the rules depicted in Figure 3 if all three edges are marked.*

If $\mathcal{T}_\ell$ is a mesh and $\mathcal{M}_\ell \subset \mathcal{T}_\ell$ is a set of marked elements, we write $\mathcal{T}_{\ell+1} = \texttt{refine}(\mathcal{T}_\ell, \mathcal{M}_\ell)$ if $\mathcal{T}_{\ell+1}$ is generated by Algorithm 2.

**Remark 3.** *For a practical realization, our formulation of Algorithm 2 does need, in fact, further decision criteria: First, the initialization step (o) needs a criterion which edges of a marked element* $T \in \mathcal{M}_\ell$ *are marked for refinement. If we want to emphasize the set* $\mathcal{M}_\ell^{(0)}$*, we will write* $\texttt{refineNVB}(\mathcal{T}_\ell, \mathcal{M}_\ell, \mathcal{M}_\ell^{(0)})$*. Moreover, the actual refinement step (iii) needs a criterion to decide which type of refinement should be used for a particular element* $T \in \mathcal{T}_\ell$ *with three marked edges. However, these criteria are not important for our results and our analysis. We therefore give only certain possible examples:*

- *If, for each element* $T \in \mathcal{M}_\ell$*, only the reference edge is marked and if all elements with three marked edges are refined by three bisections, Algorithm 2 coincides with Algorithm 1.*
- *In many implementations of newest vertex bisection, all edges of marked elements* $T \in \mathcal{M}_\ell$ *are marked in step (o). This guarantees that all edges of marked elements will be halved in the generated mesh* $\mathcal{T}_{\ell+1}$*.*
- *The numerical analysis of certain convergence results for adaptive FEM relies on the so-called* interior node property *which first appeared in* [14]*, i.e. refinement of a marked element* $T \in \mathcal{M}_\ell$ *generates a new node inside of* $T$*. For this purpose, the work* [14] *introduced the following refinement rule: For a marked element* $T \in \mathcal{M}_\ell$*, all edges are marked for refinement in step (o). For an element* $T \in \mathcal{T}_\ell$ *with three marked edges in step (iii), one uses bisec(5)-refinement in case of* $T \in \mathcal{M}_\ell$ *and bisec(3)-refinement otherwise.*
- *In* [7]*, all edges of a marked element* $T \in \mathcal{M}_\ell$ *are marked in step (o), and refinement in step (iii) is restricted to red-refinement if all edges of an element* $T \in \mathcal{T}_\ell$ *are marked, see Figure 2. This provides a variant of the popular red-green-blue refinement, see* [25, Chapter 4].

*All these cases are, in particular, covered by our analysis.*

**Remark 4.**

- *If we restrict the refinement pattern in step (iii) of Algorithm 2 to the rules of Figure 1, we write* $\mathcal{T}_{\ell+1} = \texttt{refineNVB3}(\mathcal{T}_\ell, \mathcal{M}_\ell)$*. Note that* $\mathcal{T}_{\ell+1}$ *could also be generated by two calls of Algorithm 1 via*

$$\mathcal{T}_{\ell+1/2} = \texttt{refineNVB}(\mathcal{T}_\ell, \mathcal{M}_\ell) \quad \text{and} \quad \mathcal{T}_{\ell+1} = \texttt{refineNVB}(\mathcal{T}_{\ell+1/2}, \mathcal{M}_{\ell+1/2}),$$

*where*

$$\mathcal{M}_{\ell+1/2} = \big\{T \in \mathcal{T}_{\ell+1/2} \,:\, \exists \widehat{T} \in \mathcal{M}_\ell \, \exists E \in \mathcal{M}_\ell^{(0)} \text{ with } E \subset T \subseteq \widehat{T}\big\}.$$

- *If we restrict the refinement pattern in step (iii) of Algorithm 2 to the rules of Figures 3 (left, middle) in case that all edges of an element are marked, we write* $\mathcal{T}_{\ell+1} = \texttt{refineNVBred}(\mathcal{T}_\ell, \mathcal{M}_\ell)$*.*

*Note that a call of Algorithm 2, i.e. $\mathcal{T}_{\ell+1} = \texttt{refine}(\mathcal{T}_\ell, \mathcal{M}_\ell)$, can be carried out by two calls of* `refineNVBred`*, i.e.*

$$\mathcal{T}_{\ell+1/2} = \texttt{refineNVBred}(\mathcal{T}_\ell, \mathcal{M}_\ell) \quad \text{and} \quad \mathcal{T}_{\ell+1} = \texttt{refineNVBred}(\mathcal{T}_{\ell+1/2}, \mathcal{M}_{\ell+1/2})$$

*with $\#\mathcal{M}_{\ell+1/2} \leq 2\#\mathcal{M}_\ell$.*

**2.3. Statement and discussion of main results.** We make the following general assumption: Suppose that $\mathcal{T}_0$ is a given initial mesh with *arbitrary* distribution $(E_T)_{T\in\mathcal{T}_0}$ of reference edges. Let $\mathcal{T}_\ell = \texttt{refine}(\mathcal{T}_{\ell-1}, \mathcal{M}_{\ell-1})$ for $\ell \in \mathbb{N}$ be generated inductively by Algorithm 2 with *arbitrary* sets $\mathcal{M}_j \subset \mathcal{T}_j$ of marked elements. Our first result is concerned with the mesh-closure step.

**Theorem 5** (Quasi-optimality of mesh-closure)**.** *Under the general assumption stated above, the number $\#\mathcal{T}_\ell$ of elements in the $\ell$-th mesh $\mathcal{T}_\ell$ is controlled by*

$$\#\mathcal{T}_\ell - \#\mathcal{T}_0 \leq C_1 \sum_{j=0}^{\ell-1} \#\mathcal{M}_j \quad \text{for all } \ell \in \mathbb{N}. \tag{4}$$

*The constant $C_1 > 0$ depends only on the shape regularity constant of the initial mesh $\mathcal{T}_0$, but is, in particular, independent of the sets $\mathcal{M}_j$ of marked elements chosen.*

Such a theorem has first been proved in [3, Theorem 2.4] by Binev, Dahmen, and DeVore for the usual 2D newest vertex bisection of Algorithm 1, however, under the additional assumption that the distribution $(E_T)_{T\in\mathcal{T}_0}$ of the reference edges in the initial mesh is of ***BDD-type***: For a given mesh $\mathcal{T}$, we say that the distribution $(E_T)_{T\in\mathcal{T}}$ of the reference edges is of BDD-type, if for all neighboring elements $T_1, T_2 \in \mathcal{T}$ which share an edge $E = T_1 \cap T_2$, this edge $E$ is either the reference edge of both, $T_1$ and $T_2$, or of none. In [24], Stevenson generalized this result to newest vertex bisection in $\mathbb{R}^d$ and $d \geq 2$.

For $\Omega \subseteq \mathbb{R}^2$, it is proved in [3, Lemma 2.1] that each mesh $\mathcal{T}$ admits a BDD-type choice of the reference edges. The proof is based on arguments from graph theory and does only provide an existence result. To the best of our knowledge, no efficient algorithm is known which generates such a BDD-type labeling in linear complexity. Moreover, the argument does *not* carry over to arbitrary dimension $\mathbb{R}^d$ and $d > 2$.

Our second result is concerned with the $H^1$-stability of the $L_2$-projection $\Pi_\ell : L_2(\Omega) \to \mathcal{S}^1(\mathcal{T}_\ell)$. Here,

$$\mathcal{S}^1(\mathcal{T}_\ell) := \{V_\ell : \Omega \to \mathbb{R} \ \ \mathcal{T}_\ell\text{-piecewise affine and globally continuous}\} \tag{5}$$

denotes the space of Courant finite elements of lowest-order, so-called P1-finite elements. Moreover, $\Pi_\ell$ is the $L_2$-orthogonal projection characterized by

$$\int_\Omega V_\ell \, (1 - \Pi_\ell) u \, dx = 0 \quad \text{for all } V_\ell \in \mathcal{S}^1(\mathcal{T}_\ell). \tag{6}$$

As is well-known from basic functional analysis, $\Pi_\ell$ is a continuous projection with respect to the $L_2$-norm and has operator norm 1, i.e.

$$\|\Pi_\ell u\|_{L_2(\Omega)} \leq \|u\|_{L_2(\Omega)} \quad \text{for all } u \in L_2(\Omega). \tag{7}$$

Since $\Pi_\ell$ is, in particular, well-defined on $H^1(\Omega) \supset \mathcal{S}^1(\mathcal{T}_\ell)$, it is natural to ask for $H^1$-stability

$$\|\Pi_\ell u\|_{H^1(\Omega)} \leq C_3 \, \|u\|_{H^1(\Omega)} \quad \text{for all } u \in H^1(\Omega) \tag{8}$$

with an $\ell$-independent constant $C_3 > 0$. The following theorem gives a positive answer.

**Theorem 6** ($H^1$-stability of $L_2$-projection)**.** *Under the general assumption stated above, the $L_2$-projection $\Pi_\ell$ onto $\mathcal{S}^1(\mathcal{T}_\ell)$ satisfies*

$$\|\nabla \Pi_\ell u\|_{L_2(\Omega)} \leq C_2 \, \|\nabla u\|_{L_2(\Omega)} \quad \text{for all } u \in H^1(\Omega). \tag{9}$$

*The constant $C_2 > 0$ depends only on the shape regularity constant of the initial mesh $\mathcal{T}_0$, but is independent of the sets $\mathcal{M}_j$ of marked elements chosen. In particular,* (8) *holds with $C_3 = \max\{1, C_2\}$.*

For uniform mesh-refinement, i.e. globally quasi-uniform meshes, the stability estimate (9) is well-known, see e.g. [5]: We can exploit any Clément-type quasi-interpolation operator $J_\ell$, e.g. the original operator from [9] or the Scott-Zhang projection from [19], to see

$$\begin{aligned} \|\nabla \Pi_\ell u\|_{L_2(\Omega)} &\leq \|\nabla(\Pi_\ell - J_\ell)u\|_{L_2(\Omega)} + \|\nabla J_\ell u\|_{L_2(\Omega)} \\ &\lesssim \|h_\ell^{-1}\|_{L^\infty(\Omega)} \, \|(\Pi_\ell - J_\ell)u\|_{L_2(\Omega)} + \|u\|_{L_2(\Omega)}, \end{aligned}$$

where we use the stability properties of $J_\ell$ and an inverse inequality with $h_\ell$ being the $\mathcal{T}_\ell$-piecewise constant mesh-size function $h_\ell|_T = \operatorname{diam}(T)$ for all $T \in \mathcal{T}_\ell$. Since $\Pi_\ell$ is the $L_2$-orthogonal projection, we conclude

$$\|(\Pi_\ell - J_\ell)u\|_{L(\Omega)} = \|\Pi_\ell\,(1 - J_\ell)\,u\|_{L_2(\Omega)} \leq \|(1 - J_\ell)u\|_{L_2(\Omega)} \lesssim \|h_\ell\|_{L_\infty(\Omega)} \, \|\nabla u\|_{L_2(\Omega)}.$$

Therefore, (9) holds with $C_2 \simeq \|h_\ell^{-1}\|_{L^\infty(\Omega)} \|h_\ell\|_{L^\infty(\Omega)}$, which is bounded uniformly for globally quasi-uniform meshes. We stress that we have used the stability of $\Pi_\ell$ which is only known for the non-weighted $L_2$-norm. Therefore, this proof does not cover adaptively generated locally refined meshes.

In [4], $H^1$-stability (9) is shown for meshes with a controlled local change of the elements' area. In [21], the result of [4] has been extended to $H^s$-stability $\|\Pi_\ell v\|_{H^s(\Omega)} \leq C\,\|v\|_{H^s(\Omega)}$ for all $v \in H^s(\Omega)$ and $0 \leq s \leq 1$. The conditions in [10] and [4] were finally merged in [6]. We stress that a standard interpolation argument proves that $H^1$-stability of $\Pi_\ell$ also implies $H^s$-stability for all $0 \leq s < 1$.

The first work that avoided any growth conditions for the local mesh-size or element areas, has been published by Carstensen [7]. He considered the modified NVB algorithm in 2D of Algorithm 2, where elements with three marked edges are refined by the so-called *red refinement* which splits the triangle into four similar sons, see Figure 2. Under a certain weakened BDD-assumption on the distribution $(E_T)_{T\in\mathcal{T}_0}$ of the reference edges in the initial mesh, he proved (9), see [7, Theorem 1]. We remark that [7] provides an algorithm which generates the required distribution of the reference edges for the initial mesh $\mathcal{T}_0$ in linear complexity. Moreover and with minor modifications of the analysis, similar results to those of [7] can also be obtained for the usual NVB algorithm or the general form of Algorithm 2 which is the topic of this work. We stress, however, that our analysis even removes the weakened BDD-assumption on $\mathcal{T}_0$ from [7].

## 3. Preliminaries

**3.1. Notation & abbreviations.** This short section collects some notations and definitions which are used throughout the remainder of this paper.

- We write $\mathcal{T}_\star = \texttt{refineNVB}(\mathcal{T}_\ell)$ if $\mathcal{T}_\star$ is generated by finitely many calls of Algorithm 1 from $\mathcal{T}_\ell$, i.e. there are finitely many meshes $\mathcal{T}_k, \dots, \mathcal{T}_{k+n}$ with arbitrary but fixed $n \in \mathbb{N}$ and sets of marked elements $\mathcal{M}_{k+j} \subseteq \mathcal{T}_{k+j}$ such that $\mathcal{T}_k = \mathcal{T}_\ell$, $\mathcal{T}_{k+j} = \texttt{refineNVB}(\mathcal{T}_{k+j-1}, \mathcal{M}_{k+j-1})$ for all $j = 1, \dots, n$, and finally $\mathcal{T}_\star = \mathcal{T}_{k+n}$.
- If $\mathcal{T}_\star$ is generated by finitely many calls of Algorithm 2 from $\mathcal{T}_\ell$, we write $\mathcal{T}_\star = \texttt{refine}(\mathcal{T}_\ell)$.

From now on, we fix an initial mesh $\mathcal{T}_0$ and assume that all further meshes $\mathcal{T}_\star$ are obtained from this, i.e. $\mathcal{T}_\star = \texttt{refineNVB}(\mathcal{T}_0)$, $\mathcal{T}_\star = \texttt{refineNVB3}(\mathcal{T}_0)$, $\mathcal{T}_\star = \texttt{refineNVBred}(\mathcal{T}_0)$, or $\mathcal{T}_\star = \texttt{refine}(\mathcal{T}_0)$.

- The ***set of edges*** of a mesh $\mathcal{T}_\star$ is always denoted by $\mathcal{E}_\star$.
- The ***set of nodes*** of $\mathcal{T}_\star$ is always denoted by $\mathcal{N}_\star$.
- To each element $T$, we inductively assign a ***level*** $\mathrm{gen}(T) \in \mathbb{N}_0$ as follows: For $T \in \mathcal{T}_0$, we define $\mathrm{gen}(T) := 0$. If an element $T$ is bisected into two sons $T', T''$, we define $\mathrm{gen}(T') := \mathrm{gen}(T) + 1 =: \mathrm{gen}(T'')$. If an element $T$ is red-refined into four sons $T_1, \dots, T_4$, we define $\mathrm{gen}(T_i) := \mathrm{gen}(T) + 2$. Note that one call of `refineNVB` can locally increase the level by 2 if two or three edges of $T$ are marked for refinement, see Figure 1.
- If $T \in \mathcal{T}_\star$ shares its reference edge $E_T$ with an element $T' \in \mathcal{T}_\star$ (which might have a difference reference edge), this ***reference neighbor*** is denoted by $N(T) := T'$. Otherwise, we define $N(T) := \emptyset$ as well as $N(\emptyset) := \emptyset$.
- Two elements $T_1, T_2 \in \mathcal{T}_\star$ which share an edge $E = T_1 \cap T_2 \in \mathcal{E}_\star$ are ***compatibly divisible*** if and only if the respective reference edges either satisfy $E_{T_1} = E = E_{T_2}$ or $E_{T_1} \neq E \neq E_{T_2}$, i.e. it holds $N(T_1) = T_2 \iff N(T_2) = T_1$.
- A mesh $\mathcal{T}_\star$ has the ***BDD-property*** if all pairs of elements which share an edge are compatibly divisible.
- To state the weak BDD-property of [7], we need the following definition: An element $T$ is an ***isolated element*** if $N(N(T)) \neq T$, i.e. the reference neighbor of the element $T$ has a different reference edge.
- A mesh $\mathcal{T}_\star$ has the ***weak BDD-property*** if two distinct isolated elements $T_1, T_2 \in \mathcal{T}_\star$ cannot share an edge.
- We denote by $\mathrm{chain}(\mathcal{T}_\ell, T)$ the (finite) sequence of distinct elements in $\mathcal{T}_\ell$ that are bisected by call of $\texttt{refineNVB}(\mathcal{T}_\ell, \{T\})$. Put differently, $\mathrm{chain}(\mathcal{T}_\ell, T) = (T_j)_{j=1}^n$ with $T_1 = T$, $N(T_j) = T_{j+1}$ for $j = 1, \dots, n-1$. The number $n \in \mathbb{N}$ is such that either $N(T_n) = \emptyset$ or $N(T_n)$ is already contained in $(T_j)_{j=1}^{n-1}$.
- Finally, we define $\mathcal{F}(\mathcal{T}_\star) := \big\{(T, E) \in \mathcal{T}_\star \times \mathcal{E}_\star \,:\, E \subset T\big\}$.

**3.2. Auxiliary results for NVB (Algorithm 1).** The first lemma essentially recalls some results from the literature, which are improved by some bootstrapping arguments below.

**Lemma 7.** *Assume that the initial mesh $\mathcal{T}_0$ has a BDD-type distribution of the reference edges. Then, an arbitrary refinement $\mathcal{T}_\star = \texttt{refineNVB}(\mathcal{T}_0)$ of $\mathcal{T}_0$ has the following properties:*
(i) *For each $T \in \mathcal{T}_\star$ with $N(T) \neq \emptyset$, one of the following cases holds:*

  (I) $\mathrm{gen}(T) = \mathrm{gen}(N(T))$ *and $T$, $N(T)$ are compatibly divisible, or*
  (II) $\mathrm{gen}(T) = \mathrm{gen}(N(T)) + 1$ *and $T$ is compatibly divisible with a son of $N(T)$.*

(ii) *It holds*

$$|\mathrm{gen}(T_1) - \mathrm{gen}(T_2)| \leq 1 \quad \textit{for } T_1, T_2 \in \mathcal{T}_\star \textit{ with } T_1 \cap T_2 = E \in \mathcal{E}_\star.$$

*Proof.* Statement (i) is proved in [3]. It thus only remains to prove (ii): Provided that $E$ is either the reference edge of $T_1$ or $T_2$, the claim follows directly from (i). Thus, suppose that $E$ is neither the reference edge of $T_1$ nor of $T_2$. In the uniform *bisec*(1)-refinement $\widehat{\mathcal{T}}_\star = \texttt{refineNVB}(\mathcal{T}_\star, \mathcal{T}_\star)$ of $\mathcal{T}_\star$, $E$ is then still an edge. It is even the reference edge of the two elements adjacent to it which are sons of $T_1$ and $T_2$, respectively. The level of these two elements thus coincides according to (I). Therefore, the level of $T_1$ and $T_2$ coincides. $\square$

The following lemma contains a simple observation, but it turns out to be crucial for the remainder of our analysis in the sense that it allows to reuse results from [3, 24, 16].

**Lemma 8.** *Let* $\mathcal{U}_0 \subseteq \mathcal{T}_0$ *and* $\omega := \bigcup \mathcal{U}_0 := \big\{x \in \overline{\Omega} \,:\, \exists T \in \mathcal{U}_0 \quad x \in T\big\} \subseteq \overline{\Omega}$. *Then,*

$$\mathcal{T}_\star = \texttt{refineNVB}(\mathcal{T}_0) \quad \Longrightarrow \quad \mathcal{T}_\star|_\omega = \texttt{refineNVB}(\mathcal{T}_0|_\omega),$$

*where e.g.* $\mathcal{T}_\star|_\omega := \big\{T \in \mathcal{T}_\star \,:\, T \subseteq \omega\big\}$ *denotes the restricted mesh.* $\square$

With the help of this observation, we may generalize Lemma 7 (i)–(ii) to meshes with an arbitrary distribution of reference edges. The following technical proposition is somehow the core of our improved analysis.

**Proposition 9.** *An arbitrary refinement* $\mathcal{T}_\star = \texttt{refineNVB}(\mathcal{T}_0)$ *of the initial mesh* $\mathcal{T}_0$ *has the following properties, independently of the distribution* $(E_T)_{T\in\mathcal{T}_0}$ *of the reference edges:*
(i) *There exist constants* $C_{\rm diam}, C^{\rm diam} > 0$ *which depend only on* $\mathcal{T}_0$, *such that*

$$C_{\rm diam} 2^{-{\rm gen}(T)/2} \le |T|^{1/2} \le {\rm diam}(T) \le C^{\rm diam} 2^{-{\rm gen}(T)/2} \qquad \textit{for all } T \in \mathcal{T}_\star.$$

(ii) *For* $T_1, T_2 \in \mathcal{T}_\star$ *with* $T_1 \cap T_2 = E \in \mathcal{E}_\star$ *holds* $|{\rm gen}(T_1) - {\rm gen}(T_2)| \le 2$.
(iii) *For some fixed* $T \in \mathcal{T}_\star$ *with* $N(T) \neq \emptyset$, *suppose that* ${\rm gen}(N(T)) > {\rm gen}(T)$. *Then,* $T$ *and* $N(T)$ *are compatibly divisible and* ${\rm gen}(N(T)) = {\rm gen}(T) + 1$.
(iv) *Neighboring successors* $T_1, T_2 \in \mathcal{T}_\star$ *of* $\widehat{T} \in \mathcal{T}_0$ *with the same level* ${\rm gen}(T_1) = {\rm gen}(T_2)$ *are compatibly divisible.*
(v) *Let* $\widehat{T}_1, \widehat{T}_2 \in \mathcal{T}_0$ *be compatibly divisible. Then, neighboring successors* $T_1, T_2 \in \mathcal{T}_\star$ *of* $\widehat{T}_1$ *and* $\widehat{T}_2$ *with the same level* ${\rm gen}(T_1) = {\rm gen}(T_2)$ *are compatibly divisible.*
(vi) *Suppose* ${\rm gen}(N(T)) = {\rm gen}(T)$ *and* $T$ *and* $N(T)$ *are not compatibly divisible. Then, the common edge* $E_T = T \cap N(T) \in \mathcal{E}_\star$ *is part of an edge* $E_0 \in \mathcal{E}_0$ *of the coarse mesh* $\mathcal{T}_0$.
(vii) *Let* $T \in \mathcal{T}_\star$ *and suppose that there is a finite sequence of pairwise disjoint elements* $(T_j)_{j=1}^n$ *with* $T_j \in \mathcal{T}_\star$, $T_1 = T$, $T_j = N(T_{j-1})$ *and all the* $T_j$ *have level* ${\rm gen}(T_j) = {\rm gen}(T)$. *Then, already* $n \le C_{\rm chain}$, *where the constant* $C_{\rm chain} \in \mathbb{N}$ *depends solely on* $\mathcal{T}_0$.

*Proof of Proposition 9 (i).* The statement is explicitly found in [16], but also part of the proofs in [3, 24]. $\square$

*Proof of Proposition 9 (ii).* Let $\widehat{T}_1, \widehat{T}_2 \in \mathcal{T}_0$ be the ancestors of $T_1$ resp. $T_2$ in the initial mesh $\mathcal{T}_0$. We may consider two separate cases:

**First case:** The ancestors in the initial mesh coincide $\widehat{T} := \widehat{T}_1 = \widehat{T}_2$. Then, $\mathcal{T}_0|_{\widehat{T}}$ is a BDD-mesh. Due to Lemma 8, it holds that

$$T_1, T_2 \in \mathcal{T}_\star|_{\widehat{T}} = \texttt{refineNVB}\left(\mathcal{T}_0|_{\widehat{T}}\right).$$

Hence, we can apply Lemma 7 (ii) to conclude $|{\rm gen}(T_1) - {\rm gen}(T_2)| \le 1$.

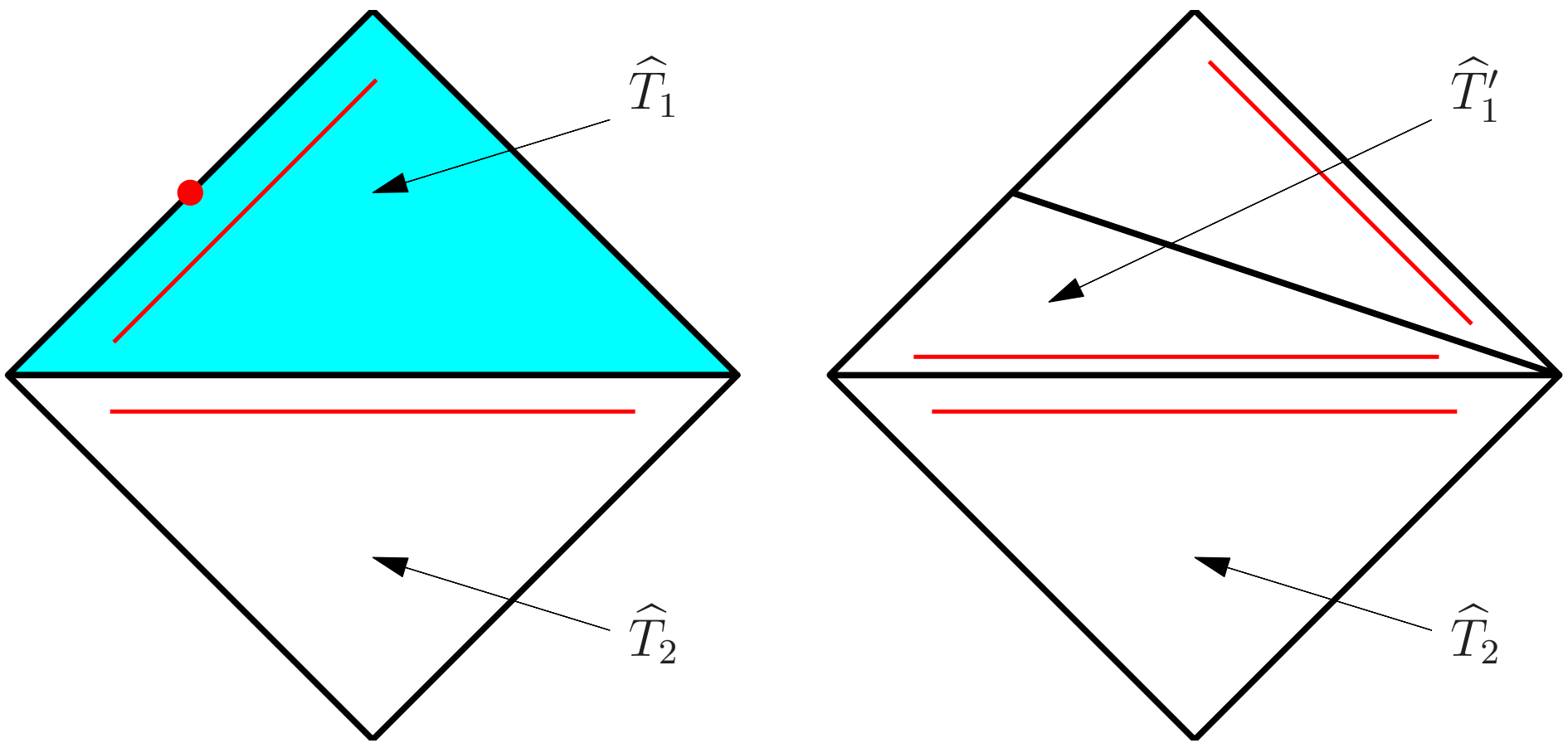

FIGURE 4. Illustration for proof of Proposition 9 (ii).

**Second case:** The elements $T_1$ and $T_2$ have different ancestors $\widehat{T}_1 \neq \widehat{T}_2$. Then, $\mathcal{E}_\star \ni E := T_1 \cap T_2 \subseteq \widehat{T}_1 \cap \widehat{T}_2 =: \widehat{E} \in \mathcal{E}_0$.

- If $\widehat{E}$ is the reference edge of both, $\widehat{T}_1$ and $\widehat{T}_2$, or if $\widehat{E}$ is neither the reference edge of $\widehat{T}_1$ nor of $\widehat{T}_2$, then $\mathcal{T}_0|_{\widehat{T}_1 \cup \widehat{T}_2}$ is a BDD-mesh. As in the first case, we conclude $|\mathrm{gen}(T_1) - \mathrm{gen}(T_2)| \leq 1$.
- The edge $\widehat{E}$ is reference edge of one of the ancestors, say $\widehat{T}_2$, but not of $\widehat{T}_1$, cf. Figure 4 (left). An arbitrary refinement $\mathcal{T}_\star$ of $\mathcal{T}_0$ that bisects either $\widehat{T}_1$ or $\widehat{T}_2$, has to bisect $\widehat{T}_1$ inevitably, which leads —at least in an intermediate state— to Figure 4 (right). Without loss of generality, we may assume $T_1 \subseteq \widehat{T}_1'$, since $T_1 = \widehat{T}_1$ would imply $T_2 = \widehat{T}_2$ and hence $\mathrm{gen}(T_1) = 0 = \mathrm{gen}(T_2)$. Next, note that the mesh $\widehat{\mathcal{T}} := \{\widehat{T}_1', \widehat{T}_2\}$ is a BDD-mesh. We introduce a new level function $\widehat{\mathrm{gen}}$ on $\widehat{\mathcal{T}}$ by $\widehat{\mathrm{gen}}(\widehat{T}_1') := 0 =: \widehat{\mathrm{gen}}(\widehat{T}_2)$. With Lemma 8, it holds that

$$T_1, T_2 \in \mathcal{T}_\star|_{\widehat{T}_1' \cup \widehat{T}_2} = \texttt{refineNVB}(\widehat{\mathcal{T}}).$$

  Lemma 7 (ii) thus yields $|\widehat{\mathrm{gen}}(T_1) - \widehat{\mathrm{gen}}(T_2)| \leq 1$. Together with $\mathrm{gen}|_{\widehat{T}_1'} = \widehat{\mathrm{gen}}|_{\widehat{T}_1'} + 1$ and $\mathrm{gen}|_{\widehat{T}_2} = \widehat{\mathrm{gen}}|_{\widehat{T}_2}$, we conclude $|\mathrm{gen}(T_1) - \mathrm{gen}(T_2)| \leq 2$.

□

*Proof of Proposition 9 (iii).* We first show that $\mathrm{gen}(N(T)) > \mathrm{gen}(T)$ ensures that $T$ and $N(T)$ are compatibly divisible. To that end, we argue by contradiction and assume that $T$ and $N(T)$ are *not* compatibly divisible. We thus start with the situation of Figure 5 (left, top) with $T = T_2$ and $N(T) = T_1$ and note that $\mathrm{gen}(T_1) > \mathrm{gen}(T_2)$ by assumption. A uniform $\mathrm{bisec}_3$-refinement of $\mathcal{T}_\ell$ leads to the situation of Figure 5 (right, top). Marking the highlighted element leads us to the situation of Figure 5 (left, bottom). Now, we proceed as follows: outside of the plotted scope, we refine the mesh by a bisec(3)-refinement, which can be carried out by two bisec(1)-steps. Additionally, we mark the highlighted element for refinement. Doing so, we finally end up with Figure 5 (right, bottom). For the elements $T_1'$ and $T_2'$ holds $\mathrm{gen}(T_1') = \mathrm{gen}(T_1) + 5$ and $\mathrm{gen}(T_2') = \mathrm{gen}(T_2) + 3$. Together with the assumption

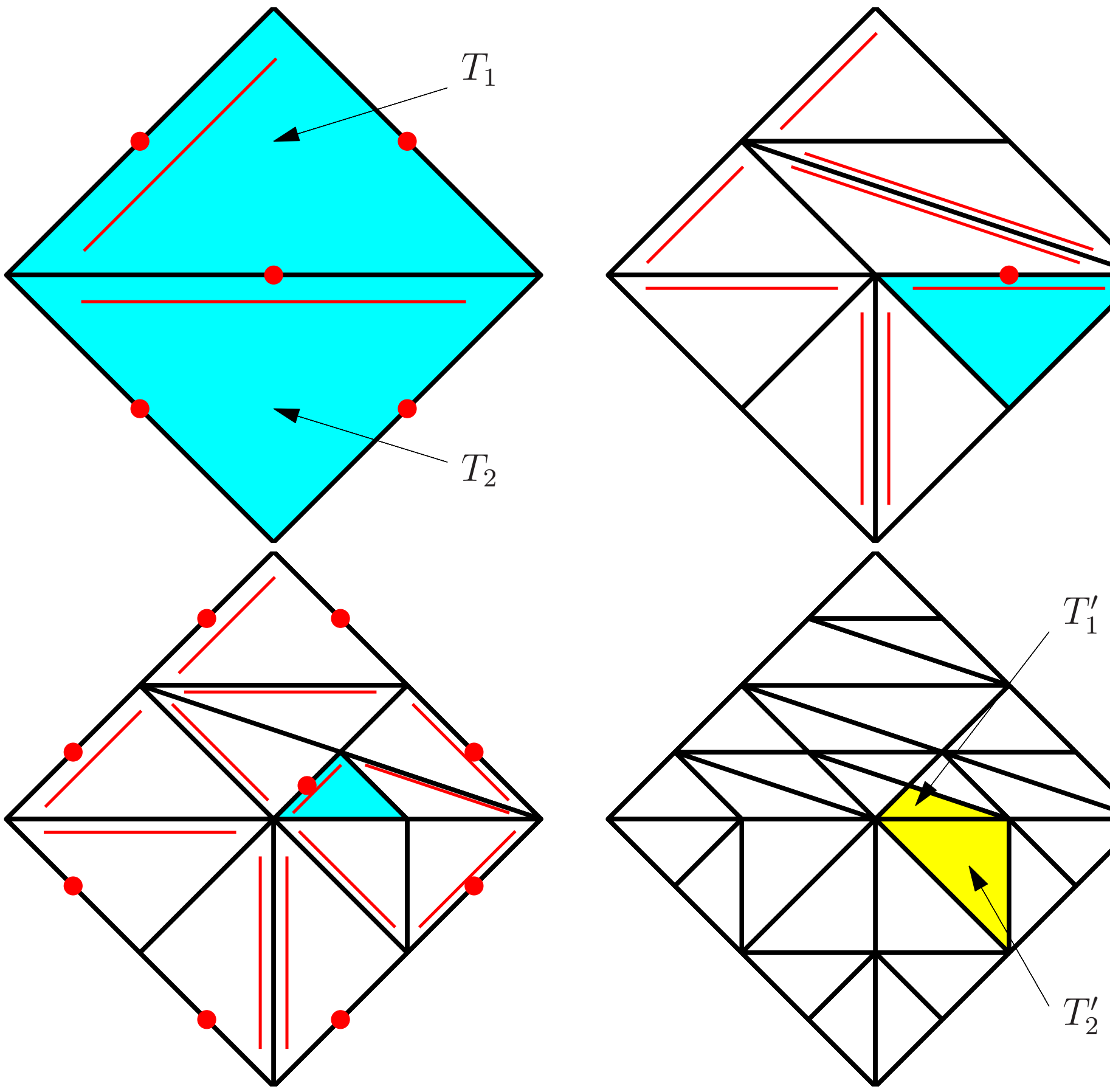

FIGURE 5. Illustriation of proof of Proposition 9 (iii).

$\mathrm{gen}(T_1) > \mathrm{gen}(T_2)$, we are thus led to

$$\mathrm{gen}(T'_1) = \mathrm{gen}(T_1) + 5 > \mathrm{gen}(T_2) + 5 = \mathrm{gen}(T'_2) + 2.$$

This contradicts point (ii) of Proposition 9 and thus proves that $T$ and $N(T)$ need to be compatibly divisible.

Second, it remains to show $\mathrm{gen}(N(T)) \leq \mathrm{gen}(T) + 1$. We again argue by contradiction: From the last step, we know that $T$ and $N(T)$ are compatibly divisible. This setting is visualized in Figure 6, where again $T = T_2$ and $N(T) = T_1$. Some local refinements lead to Figure 6 (right, bottom). There, $\mathrm{gen}(T'_1) = \mathrm{gen}(T_1) + 4$ and $\mathrm{gen}(T'_2) = \mathrm{gen}(T_2) + 3$. Therefore, $\mathrm{gen}(T_1) > \mathrm{gen}(T_2) + 1$ would yield

$$\mathrm{gen}(T'_1) = \mathrm{gen}(T_1) + 4 > \mathrm{gen}(T_2) + 5 = \mathrm{gen}(T'_2) + 2$$

and hence a contraction to point (ii) of Proposition 9. Consequently, it holds that $\mathrm{gen}(N(T)) \leq \mathrm{gen}(T) + 1$. □

*Proof of Proposition 9 (iv)–(vii).* (iv) Let $\widehat{T} \in \mathcal{T}_0$ be the ancestor of $T \in \mathcal{T}_\star$ and note that $\mathcal{T}_0|_{\widehat{T}}$ is a BDD-mesh with $T \in \mathcal{T}_\star|_{\widehat{T}} = \texttt{refineNVB}(\mathcal{T}_0|_{\widehat{T}})$. The statement thus follows with Lemma 7.

(v) This follows analogously to (iv) since $\mathcal{T}_0|_{\widehat{T}_1 \cup \widehat{T}_2}$ is a BDD-mesh.

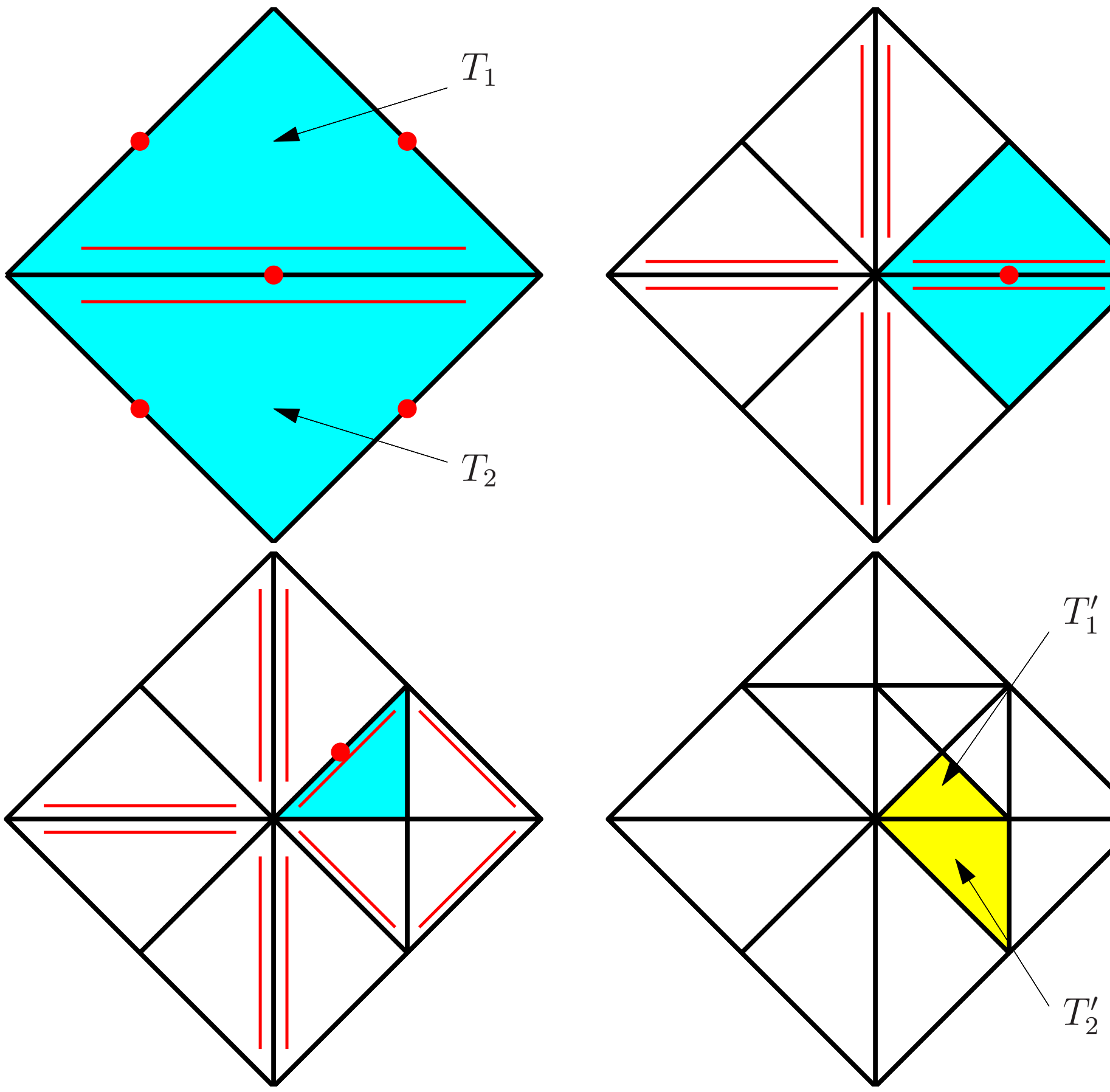

FIGURE 6. Illustration of proof of Proposition 9 (iii).

(vi) We show the contraposition: If $E$ is not a subset of a coarse edge, it lies inside of a coarse element $\widehat{T} \in \mathcal{T}_0$. By use of $\mathrm{gen}(T) = \mathrm{gen}(N(T))$, (iv) predicts that $T$ and $N(T)$ are compatibly divisible.

(vii) Without loss of generality, we assume $n \geq 4$. By assumption, the pairs $(T_i, T_{i+1})$ for $i \in \{1, \dots, n-2\}$ cannot be compatibly divisible. Furthermore, we claim that the reference edges of $T_i$ for $i \in \{1, \dots, n-3\}$ contain at least one node from $\mathcal{T}_0$: If $T_k$ was an element with a reference edge not sharing a node in $\mathcal{T}_0$, we immediately see that the other two edges of $T_{k+1}$ would not be parts of coarse edges. But then, $T_{k+1}$ and $T_{k+2}$ would have the same ancestor from $\mathcal{T}_0$. Hence, $T_{k+1}$ and $T_{k+2}$ would be compatibly divisible according to (vi). We thus conclude that the reference edges of the elements $T_i$ for $i \in \{1, \dots, n-3\}$ contain at least one node of the coarse mesh $\mathcal{T}_0$. The newest-vertex bisection algorithm yields uniformly shape-regular meshes. Therefore, the number of elements touching a node is bounded uniformly. This leads us to

$$n - 3 \lesssim \#\mathcal{N}_0$$

and hence to $n \leq C$. The definition $C_{\mathrm{chain}} := \lceil C \rceil$ concludes the proof. $\square$

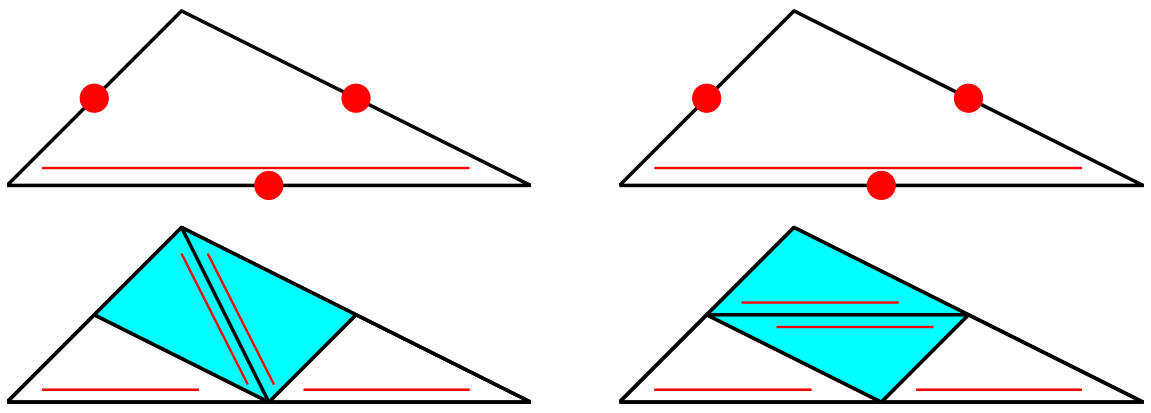

FIGURE 7. An element $T$ with 3 marked edges is refined with bisec(3) when using NVB (left). When using MNVB of Section 3.4, $T$ is either refined by bisec(3) or red-refined into 4 congruent sons (right). The so-called bisec(3)-sons as well as the red-sons are highlighted.

**3.3. Auxiliary results for MNVB with bisec(3)-refinement only (Algorithm 2).** In this section, we consider `refineNVB3`. As stated in Remark 4, we consider Algorithm 2 under the following restrictions: In the refinement step (iii), only bisec(3)-refinement is used if all edges of an element $T \in \mathcal{T}_\ell$ are marked, i.e. we stick with the usual refinement pattern of NVB from Figure 1. Compared to `refineNVB` of Algorithm 1, the only difference is that in step (o) not only reference edges are marked for refinement.

Analogously to the notation introduced before, we write $\mathcal{T}_\star :=$ `refineNVB3`$(\mathcal{T}_\ell)$. The essential observation from Remark 4 is the following: Each mesh $\mathcal{T}_\star =$ `refineNVB3`$(\mathcal{T}_0)$ also satisfies $\mathcal{T}_\star =$ `refineNVB`$(\mathcal{T}_0)$, since each step of bisec(3)-NVB can be done by two successive steps of the usual NVB algorithm. In particular, we obtain the following result:

**Proposition 10.** *Proposition 9 also holds if* `refineNVB3` *is used instead of* `refineNVB`*.* $\square$

**3.4. Auxiliary results for MNVB with bisec(3)- and red-refinement (Algorithm 2).** In this section, we consider Algorithm 2 under the following restriction on the refinement step (iii): If all edges of an element $T \in \mathcal{T}_\ell$ are marked, the element is either refined by bisec(3)-refinement or by red-refinement, i.e. we still exclude the bisec(5)-refinement from the refinement pattern of Figure 3. The two sons of a red-refined triangle $T$, that do not intersect the reference edge $E_T$ of $T$ in more than one point, are called (lower and upper) ***red-sons*** of $T$, see Figure 7 (right). Note that red-sons are the only elements that cannot be created by the usual newest vertex bisection algorithm.

We will need a special feature of the set of marked edges $\mathcal{M}_\ell^{(k)}$ that is computed in Algorithm 2.

**Lemma 11.** *The set $\mathcal{M}_\ell^{(k)}$ computed by Algorithm 2 is, with respect to cardinality and even set inclusion, the smallest set with the properties*

$$\mathcal{M}_\ell^{(0)} \subset \mathcal{M}_\ell^{(k)} \tag{10}$$
$$\text{if an edge } E \text{ of an element } T \text{ fulfills } E \in \mathcal{M}_\ell^{(k)}, \text{ then } E_T \in \mathcal{M}_\ell^{(k)}.$$

*In particular, $\mathcal{M}_\ell^{(k)}$ is the unique set of minimal cardinality that fulfills* (10). $\square$

**Definition 12.** *We call the meshes $\mathcal{T} = (\mathcal{T}, (E_T)_{T\in\mathcal{T}})$ and $\widetilde{\mathcal{T}} = (\widetilde{\mathcal{T}}, (E_T)_{T\in\widetilde{\mathcal{T}}})$ corresponding meshes, written $\mathcal{T} \simeq \widetilde{\mathcal{T}}$, if there exists a bijective function* $\mathrm{corr} : \mathcal{F}(\mathcal{T}) \to \mathcal{F}(\widetilde{\mathcal{T}})$ *with the* *following properties: For* $(T,E),(T',E')\in\mathcal{F}(\mathcal{T})$ *with*

$$(\widetilde{T},\widetilde{E}) := \operatorname{corr}(T,E) \quad \text{and} \quad (\widetilde{T}',\widetilde{E}') := \operatorname{corr}(T',E'),$$

*holds*

(i) $|T|\simeq|\widetilde{T}|$ *as well as* $\operatorname{gen}(T)=\operatorname{gen}(\widetilde{T})$.

(ii) *It holds* $\widetilde{T}\cap\widetilde{T}'=\widetilde{E}\in\widetilde{\mathcal{E}}$ *if and only if* $T\cap T'=E\in\mathcal{E}$*, i.e. neighboring elements are mapped to neighboring elements.*

(iii) *It holds* $\widetilde{E}=E_{\widetilde{T}}$ *if and only if* $E=E_T$*, i.e. reference edges are mapped to reference edges.*

(iv) *It holds* $E=E'=E_T$ *and* $T'=N(T)$ *if and only if* $\widetilde{E}=\widetilde{E}'=E_{\widetilde{T}}$ *and* $\widetilde{T}'=N(\widetilde{T})$.

(v) $\widetilde{T}$ *and* $\widetilde{T}'$ *are compatibly divisible if and only if* $T$ *and* $T'$ *are compatibly divisible.*

(vi) $\widetilde{T}$ *and* $\widetilde{T}'$ *are neighboring successors of an ancestor* $\widehat{T}\in\mathcal{T}_0$ *if and only if* $T$ *and* $T'$ *are so.*

(vii) $T'=T$ *and* $E'=E_T$ *implies* $E_{\widetilde{T}'}=\widetilde{E}'=E_{\widetilde{T}}$*, i.e.* $\widetilde{E}'$ *is the reference edge of both* $\widetilde{T}$ *and* $\widetilde{T}'$*. Furthermore,* $\widetilde{T}'=\widetilde{T}$ *and* $\widetilde{E}'=E_{\widetilde{T}}$ *implies* $E_{T'}=E'=E_T$*, i.e.* $E'$ *is the reference edge of both* $T$ *and* $T'$*.*

**Lemma 13.** *Suppose that* $\mathcal{T}$ *and* $\widetilde{\mathcal{T}}$ *are corresponding meshes* $\mathcal{T}\simeq\widetilde{\mathcal{T}}$ *and* $T\in\mathcal{T}$*. Then,*

- $\#\mathcal{T}=\#\widetilde{\mathcal{T}}$*, and*
- $\#\operatorname{corr}(T)\le 2$ *for* $\operatorname{corr}(T):=\{\widetilde{T}\in\widetilde{\mathcal{T}}\,:\,(\widetilde{T},\widetilde{E})=\operatorname{corr}(T,E)$ *for some* $E\in\mathcal{E}_\ell$ *with* $E\subset T\}$.

*Proof.* Since $\operatorname{corr}:\mathcal{F}(\mathcal{T})\to\mathcal{F}(\widetilde{\mathcal{T}})$ is bijective between two finite sets, we have $\#\mathcal{F}(\mathcal{T})=\#\mathcal{F}(\widetilde{\mathcal{T}})$. On the other hand, $\#\mathcal{F}(\mathcal{T})=3\times\#\mathcal{T}$ and $\#\mathcal{F}(\widetilde{\mathcal{T}})=3\times\#\widetilde{\mathcal{T}}$, whence $\#\mathcal{T}=\#\widetilde{\mathcal{T}}$.

To show the second statement, suppose that $(\widetilde{T},E_{\widetilde{T}})=\operatorname{corr}(T,E_T)$. Now choose an edge $E$ of $T$ that is not the reference edge and consider $(\widetilde{T}',\widetilde{E}')=\operatorname{corr}(T,E)$ with $\widetilde{T}'\neq\widetilde{T}$. From Definition 12, (vii), we know that $E_{\widetilde{T}}$ is the reference edge of $\widetilde{T}'$. Hence, $\widetilde{T}\cap\widetilde{T}'=E_{\widetilde{T}}$. ☐

**Lemma 14.** *Denote by* $(\mathcal{T}_\ell)_{\ell\in\mathbb{N}_0}$ *a sequence of meshes which is inductively generated by* `refineNVBred` *from a coarse mesh* $\mathcal{T}_0$*, i.e.*

$$\mathcal{T}_{\ell+1}=\texttt{refineNVBred}(\mathcal{T}_\ell,\mathcal{M}_\ell,\mathcal{M}_\ell^{(0)}) \quad \text{for } \ell\in\mathbb{N}_0$$

*with certain marked elements* $\mathcal{M}_\ell\subseteq\mathcal{T}_\ell$ *and corresponding marked edges* $\mathcal{M}_\ell^{(0)}\subseteq\{E\in\mathcal{E}_\ell : E\subset T$ *for some* $T\in\mathcal{M}_\ell\}$*. Then, there is a sequence of corresponding meshes* $(\widetilde{\mathcal{T}}_\ell)_{\ell\in\mathbb{N}_0}$*, i.e.* $\mathcal{T}_\ell\simeq\widetilde{\mathcal{T}}_\ell$*, which is inductively generated by* `refineNVB3` *, i.e.:*

$$\widetilde{\mathcal{T}}_{\ell+1}=\texttt{refineNVB3}(\widetilde{\mathcal{T}}_\ell,\widetilde{\mathcal{M}}_\ell,\widetilde{\mathcal{M}}_\ell^{(0)}) \quad \text{for } \ell\in\mathbb{N}_0$$

*such that* $\#\widetilde{\mathcal{M}}_\ell\le 2\#\mathcal{M}_\ell$ *for* $\ell\in\mathbb{N}_0$*. Furthermore, if* $T\in\mathcal{T}_\ell$ *is a red-son, the map* $\operatorname{corr}_\ell:\mathcal{F}(\mathcal{T}_\ell)\to\mathcal{F}(\widetilde{\mathcal{T}}_\ell)$ *behaves as in Figure 8.*

*Proof.* We argue by induction on $\ell$. For $\ell=0$, we define $\widetilde{\mathcal{T}}_0:=\mathcal{T}_0$ and let $\operatorname{corr}_0$ be the identity. Now, suppose that $\widetilde{\mathcal{T}}_\ell$ is corresponding to $\mathcal{T}_\ell$ and that

$$\mathcal{T}_{\ell+1}=\texttt{refineNVBred}(\mathcal{T}_\ell,\mathcal{M}_\ell,\mathcal{M}_\ell^{(0)}). \tag{11}$$

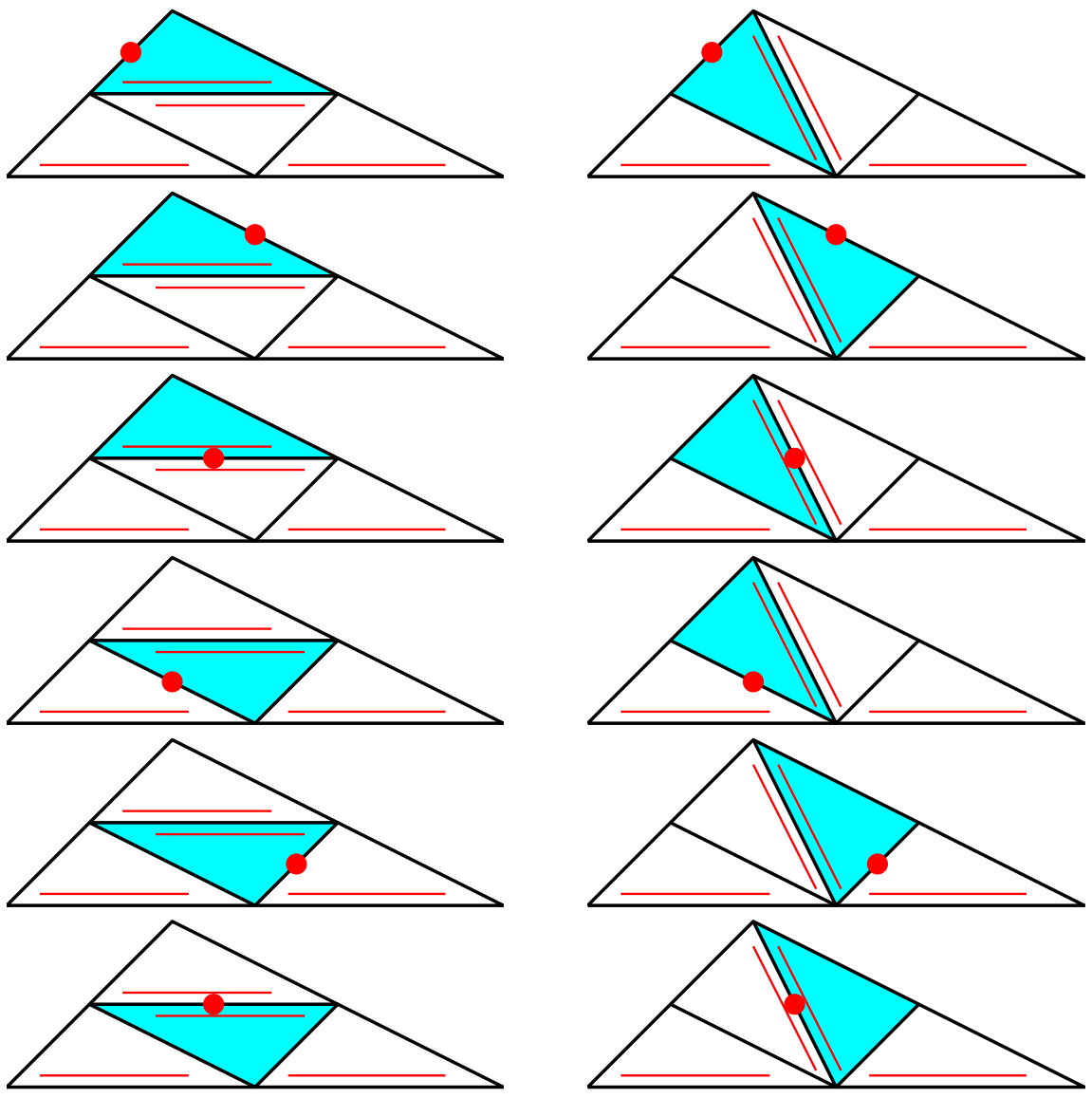

FIGURE 8. Illustration of proof of Lemma 14: Some examples of the mapping of the function $\mathrm{corr}_{\ell+1}$ on elements in $\mathcal{T}_{\ell+1} \backslash \widetilde{\mathcal{T}}_{\ell+1}$.

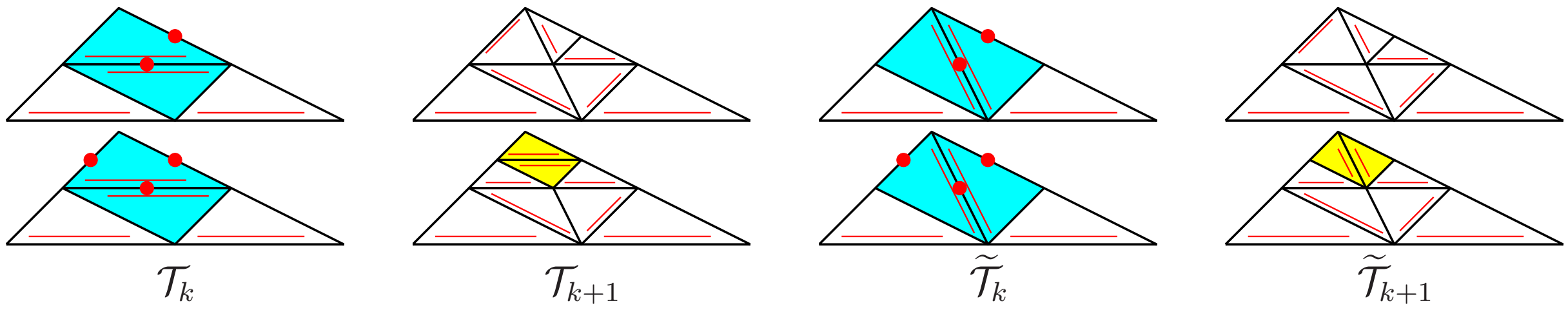

FIGURE 9. Illustration of proof of Lemma 14: In the first two columns of the first row, we see the refinement if the two red-sons are marked for refinement; one of them with an additionally marked edge. In the third and fourth column of the first row, the refinement of the two corresponding bisec(3)-sons with the same additionally marked edge is shown.

In the following, we construct a new mesh $\widetilde{\mathcal{T}}_{\ell+1}$ and a function $\mathrm{corr}_{\ell+1}$ so that $\widetilde{\mathcal{T}}_{\ell+1} \simeq \mathcal{T}_{\ell+1}$. To that end, we define the following sets:

$$\mathcal{M}_\ell^{\mathrm{corr}} := \big\{(T,E) \in \mathcal{M}_\ell \times \mathcal{M}_\ell^{(0)} \,:\, E \subseteq T\big\} \subseteq \mathcal{F}(\mathcal{T}_\ell),$$
$$\widetilde{\mathcal{M}}_\ell := \big\{\widetilde{T} \in \widetilde{\mathcal{T}}_\ell \,:\, (\widetilde{T}, \widetilde{E}) \in \mathrm{corr}_\ell(\mathcal{M}_\ell^{\mathrm{corr}})\big\},$$
$$\widetilde{\mathcal{M}}_\ell^{(0)} := \big\{\widetilde{E} \in \widetilde{\mathcal{E}}_\ell \,:\, (\widetilde{T}, \widetilde{E}) \in \mathrm{corr}_\ell(\mathcal{M}_\ell^{\mathrm{corr}})\big\}.$$

From Lemma 13, we conclude that $\#\widetilde{\mathcal{M}}_\ell \le 2\#\mathcal{M}_\ell$. Moreover, for every $\widetilde{E} \in \widetilde{\mathcal{M}}_\ell^{(0)}$ there is a $\widetilde{T} \in \widetilde{\mathcal{T}}_\ell$ such that $\widetilde{E} \subset \widetilde{T}$ and $(\widetilde{T}, \widetilde{E}) \in \widetilde{\mathcal{M}}_\ell^{\mathrm{corr}}$. We may therefore define

$$\widetilde{\mathcal{T}}_{\ell+1} = \texttt{refineNVB3}(\widetilde{\mathcal{T}}_\ell, \widetilde{\mathcal{M}}_\ell, \widetilde{\mathcal{M}}_\ell^{(0)}). \tag{12}$$

Now, consider the set of marked edges $\mathcal{M}_\ell^{(k)}$ implicitly generated by (11) and define the sets

$$\mathcal{R}_\ell^{\rm corr} := \big\{(T,E) \in \mathcal{T}_\ell \times \mathcal{M}_\ell^{(k)} \,:\, E \subset T\big\} \subseteq \mathcal{F}(\mathcal{T}_\ell),$$
$$\widetilde{\mathcal{X}}_\ell := \big\{\widetilde{E} \in \widetilde{\mathcal{E}}_\ell \,:\, (\widetilde{T},\widetilde{E}) \in {\rm corr}_\ell(\mathcal{R}_\ell^{\rm corr})\big\}.$$

Additionally, let $\widetilde{\mathcal{M}}_\ell^{(\widetilde{k})}$ denote the set of marked edges implicitly generated by (12). In the following, we employ Lemma 11 to prove $\widetilde{\mathcal{M}}_\ell^{(\widetilde{k})} = \widetilde{\mathcal{X}}_\ell$:

- First, we show $\widetilde{\mathcal{M}}_\ell^{(0)} \subset \widetilde{\mathcal{X}}_\ell$: Let $\widetilde{E} \in \widetilde{\mathcal{M}}_\ell^{(0)}$ and choose $\widetilde{T} \in \widetilde{\mathcal{T}}_\ell$ and $(T,E) \in \mathcal{M}_\ell^{\rm corr}$ with $(\widetilde{T},\widetilde{E}) = {\rm corr}_\ell(T,E)$. From $\mathcal{M}_\ell^{\rm corr} \subseteq \mathcal{R}_\ell^{\rm corr}$, we conclude $\widetilde{E} \in \widetilde{\mathcal{X}}_\ell$.
- Second, we show that $\widetilde{\mathcal{M}}_\ell^{(\widetilde{k})} \subseteq \widetilde{\mathcal{X}}_\ell$. To that end, we show that if an element $\widetilde{T} \in \widetilde{\mathcal{T}}_\ell$ has an edge in $\widetilde{\mathcal{X}}_\ell$, then also it's reference edge is in $\widetilde{\mathcal{X}}_\ell$: Let $\widetilde{T} \in \widetilde{\mathcal{T}}_\ell$ and $\widetilde{E} \in \widetilde{\mathcal{X}}_\ell$ with $\widetilde{E} \subset \widetilde{T}$. Without loss of generality, we may assume $\widetilde{E} \neq E_{\widetilde{T}}$. By definition of $\widetilde{\mathcal{X}}_\ell$, there are $\widetilde{T}' \in \widetilde{\mathcal{T}}_\ell$ and $(T',E') \in \mathcal{R}_\ell^{\rm corr}$ with $(\widetilde{T}',\widetilde{E}) = {\rm corr}_\ell(T',E')$ and, in particular, $\widetilde{E} \subset \widetilde{T}'$. We now have to distinguish two cases, where we employ different properties of the function ${\rm corr}_\ell$ stated in Definition 12.
  - (I) $\widetilde{T}' = \widetilde{T}$: By definition of $\mathcal{R}_\ell^{\rm corr}$, $(T',E') \in \mathcal{R}_\ell^{\rm corr}$ yields $(T',E_{T'}) \in \mathcal{R}_\ell^{\rm corr}$. We define $(\widetilde{T}'',\widetilde{E}'') = {\rm corr}_\ell(T',E_{T'})$ and observe $\widetilde{E}'' \in \widetilde{\mathcal{X}}_\ell$. By (iii) and (vii), we see $E_{\widetilde{T}'} = E_{\widetilde{T}''} = \widetilde{E}''$. This yields $E_{\widetilde{T}} = E_{\widetilde{T}'} \in \widetilde{\mathcal{X}}_\ell$.
  - (II) $\widetilde{T}' \neq \widetilde{T}$ and hence $\widetilde{E} = \widetilde{T} \cap \widetilde{T}'$: Let $(T,E) \in \mathcal{F}(\mathcal{T}_\ell)$ with $(\widetilde{T},\widetilde{E}) = {\rm corr}_\ell(T,E)$. By (ii), we see $T \cap T' = E = E'$. Hence $(T,E) \in \mathcal{R}_\ell^{\rm corr}$ yields $(T,E_T) \in \mathcal{R}_\ell^{\rm corr}$. Now we proceed as before. We define $(\widetilde{T}'',\widetilde{E}'') = {\rm corr}_\ell(T,E_T)$ and observe $\widetilde{E}'' \in \widetilde{\mathcal{X}}_\ell$. By (iii) and (vii), we see $E_{\widetilde{T}} = E_{\widetilde{T}''} = \widetilde{E}''$. This yields $E_{\widetilde{T}} \in \widetilde{\mathcal{X}}_\ell$.
- Finally, we need to show that $\widetilde{\mathcal{X}}_\ell$ is the smallest set with the properties (10): To that end, we define the sets

  $$\widetilde{\mathcal{R}}_\ell^{\rm corr} := \big\{(\widetilde{T},\widetilde{E}) \in \widetilde{\mathcal{T}}_\ell \times \widetilde{\mathcal{M}}_\ell^{(\widetilde{k})} \,:\, \widetilde{E} \subset \widetilde{T}\big\},$$
  $$\mathcal{X}_\ell := \big\{E \in \mathcal{E}_\ell \,:\, (T,E) \in {\rm corr}_\ell^{-1}(\widetilde{\mathcal{R}}_\ell^{\rm corr})\big\}.$$

  Since $\widetilde{\mathcal{M}}_\ell^{(\widetilde{k})}$ satisfies (10) with respect to $\widetilde{\mathcal{T}}_\ell$ and since ${\rm corr}_\ell^{-1}$ preserves the same properties as ${\rm corr}_\ell$ does, we may argue as before to see that the set $\mathcal{X}_\ell \subseteq \mathcal{M}_\ell^{(k)}$ satisfies (10) as well. Moreover, one can prove that a strict inclusion $\widetilde{\mathcal{M}}_\ell^{(\widetilde{k})} \subsetneqq \widetilde{\mathcal{X}}_\ell$ would imply the strict inclusion $\mathcal{X}_\ell \subsetneqq \mathcal{M}_\ell^{(k)}$ and hence contradict Lemma 11.

So far, Lemma 11 applies and proves equality $\widetilde{\mathcal{X}}_\ell = \widetilde{\mathcal{M}}_\ell^{(\widetilde{k})}$. Consequently, this yields that one obtains the same mesh $\widetilde{\mathcal{T}}_{\ell+1}$ independently of whether one transfers the initially marked elements/edges via ${\rm corr}_\ell$ and uses `refineNVB3` to compute the mesh-closure with respect to $\widetilde{\mathcal{T}}_\ell$, or one computes the mesh-closure via `refine` with respect to $\mathcal{T}_\ell$ and transfers the resulting marked elements/edges via ${\rm corr}_\ell$.

It now only remains to define the corresponding function ${\rm corr}_{\ell+1}$ between $\mathcal{T}_{\ell+1}$ and $\widetilde{\mathcal{T}}_{\ell+1}$: Let $T \in \mathcal{T}_\ell$. According to Lemma 13, it holds that $\#{\rm corr}_\ell(T) \in \{1,2\}$.

- Suppose that ${\rm corr}_\ell(T) = \{\widetilde{T}\}$. Since $\widetilde{\mathcal{M}}_\ell^{(\widetilde{k})} = \widetilde{\mathcal{X}}_\ell$, $T$ and $\widetilde{T}$ have the same number of marked edges. Note that an element $T \in \mathcal{T}_\ell$ can be refined by either one, two, or three bisections or by red-refinement, whereas an element $\widetilde{T} \in \widetilde{\mathcal{T}}_\ell$ can be refined by either one, two or three bisections only. If $T$ is refined by bisections only, the local definition of $\mathrm{corr}_{\ell+1}$ is obvious. If $T$ is red-refined, the local definition of $\mathrm{corr}_{\ell+1}$ is given in Figure 8.

- If $\mathrm{corr}_\ell(T) = \{\widetilde{T}_1, \widetilde{T}_2\}$, then —according to the induction hypothesis— $T$ is a red-son. In this case, $\mathrm{corr}_\ell$ behaves as in Figure 8. Independently of the marked edges of $T$ and $N(T)$, any refinement easily yields a new definition of $\mathrm{corr}_{\ell+1}$, see Figure 9.

Clearly, the new function $\mathrm{corr}_{\ell+1}$ fulfills all the properties of Definition 12, which can be easily verified by the Figures 8 and 9. ☐

From Lemma 14 and Definition 12, (i), we obtain the following result.

**Proposition 15.** *Claims (i) and (ii) of Proposition 9 hold also for the modified NVB algorithm* `refineNVBred`. ☐

**Remark 16.** *Our following analysis only requires that (i) & (ii) of Proposition 9 also hold for* `refineNVB3`. *However, by use of Lemma 14, one can even transfer (iii)–(vi) from* `refineNVB` *to* `refineNVBred`.

**3.5. Auxiliary results for MNVB (Algorithm 2).** Finally, we note that the Claims (i) and (ii) of Proposition 9 are satisfied for Algorithm 2 without any restrictions on step (iii). This is due to the fact that a call of `refine` without any restrictions can be carried out by two calls of `refineNVBred`.

**Proposition 17.** *Claims (i) and (ii) of Proposition 9 hold for Algorithm 2.* ☐

## 4. Proof of Theorem 5 (Quasi-Optimality of Mesh-Closure)

The following auxiliary lemma, as well as the proof of Theorem 5, follow the lines of Stevenson [24], but with certain important modifications.

**Lemma 18.** *Let $\mathcal{T}_0$ be an initial mesh with arbitrary distribution $(E_T)_{T\in\mathcal{T}_0}$ of reference edges and $\mathcal{T}_\ell =$* `refineNVB`*$(\mathcal{T}_0)$ be an arbitrary refinement of $\mathcal{T}_0$. Then, for arbitrary $T \in \mathcal{T}_\ell$ and all $T'$ that are created by the call of* `refineNVB`*$(\mathcal{T}_\ell, \{T\})$, it holds that*

$$\mathrm{gen}(T') \le \mathrm{gen}(T) + 2 \tag{13}$$

*and that*

$$\mathrm{dist}(T, T') \le C_{\mathrm{dist}} C^{\mathrm{diam}} \sum_{k=\mathrm{gen}(T')}^{\mathrm{gen}(T)+2} 2^{-k/2} \le \frac{C_{\mathrm{dist}} C^{\mathrm{diam}}}{1 - 2^{-1/2}} 2^{-\mathrm{gen}(T')/2}. \tag{14}$$

*The constant $C_{\mathrm{dist}} := 2C_{\mathrm{chain}} + 2^{3/2}$ depends solely on the initial mesh $\mathcal{T}_0$.*

*Proof.* We first prove (13): For all elements $\widetilde{T}$ in $\mathrm{chain}(\mathcal{T}_\ell, T)$ with level $\mathrm{gen}(\widetilde{T}) \le \mathrm{gen}(T)$, statement (13) follows from the fact that refinement of an elements increases the level function locally by at most 2, see Figure 1. If $\widetilde{T}$ is an element of $\mathrm{chain}(\mathcal{T}_\ell, T)$ with level $\mathrm{gen}(\widetilde{T}) > \mathrm{gen}(T)$, Proposition 9 (iii) predicts that $\widetilde{T}$ is the last element in the chain and $\mathrm{gen}(\widetilde{T}) = \mathrm{gen}(T)+1$. Therefore, it is only bisected once, and so it's sons have level $\mathrm{gen}(T)+2$.

We now turn to the proof of (14). The second estimate is a simple application of the geometric series, and it only remains to prove the first: Let $\widetilde{T} \in \mathrm{chain}(\mathcal{T}_\ell, T)$ be the father of the element $T'$. If $c_\ell(T) := \#\mathrm{chain}(\mathcal{T}_\ell, T) \le 2$, it holds $\mathrm{dist}(T, T') = 0$ so that the estimate holds trivially. We may therefore additionally assume $c_\ell(T) > 2$, and argue by induction on $\mathrm{gen}(T)$:

First, let $\mathrm{gen}(T) = 0$ and $\mathrm{chain}(\mathcal{T}_\ell, T) = (T_j)_{j=1}^{c_\ell(T)}$ be the chain of $T$. According to (13), it holds $\mathrm{gen}(T') \in \{1, 2\}$. Note that Proposition 9 (vii) and (ii) predict $c_\ell(T) \le C_{\mathrm{chain}} + 1$. With the triangle inequality and Lemma 7 with $\mathrm{gen}(T_k) = 0$, we obtain

$$\mathrm{dist}(T, T') \le \sum_{k=2}^{c_\ell(T)-1} \mathrm{diam}(T_k) \le \sum_{k=2}^{c_\ell(T)-1} C^{\mathrm{diam}} \le 2C_{\mathrm{chain}} C^{\mathrm{diam}} 2^{-1}.$$

Using $2^{-1} \le \sum_{k=\mathrm{gen}(T')}^{\mathrm{gen}(T)+2} 2^{-k/2}$ and $2C_{\mathrm{chain}} \le C_{\mathrm{dist}}$, we conclude

$$\mathrm{dist}(T, T') \le C_{\mathrm{dist}} C^{\mathrm{diam}} \sum_{k=\mathrm{gen}(T')}^{\mathrm{gen}(T)+2} 2^{-k/2}.$$

This proves (14) for $\mathrm{gen}(T) = 0$.

According to the induction hypothesis, we may assume that (14) holds for all $T$ with $\mathrm{gen}(T) \le n$. Let $T \in \mathcal{T}_\ell$ with $\mathrm{gen}(T) = n + 1$. We choose the first $J$ elements $(T_j)_{j=1}^J$ of $\mathrm{chain}(\mathcal{T}_\ell, T)$ such that $\mathrm{gen}(T_j) = \mathrm{gen}(T)$ for $j = 2, \ldots, J$ and $\mathrm{gen}(T_{J+1}) \ne \mathrm{gen}(T)$ in case of $c_\ell(T) > J$. Proposition 9 (vi) already yields $J \le C_{\mathrm{chain}}$. We need to consider three cases:

(A) $\widetilde{T} \in (T_j)_{j=1}^J$,

(B) $\widetilde{T} \in \mathrm{chain}(\mathcal{T}_\ell, T) \setminus (T_j)_{j=1}^J$, $c_\ell(T) > J$, and $\mathrm{gen}(T_{J+1}) > \mathrm{gen}(T)$,

(C) $\widetilde{T} \in \mathrm{chain}(\mathcal{T}_\ell, T) \setminus (T_j)_{j=1}^J$, $c_\ell(T) > J$, and $\mathrm{gen}(T_{J+1}) < \mathrm{gen}(T)$.

In case (A), we have $\mathrm{gen}(T') \in \{\mathrm{gen}(T) + 1, \mathrm{gen}(T) + 2\}$ and $\mathrm{dist}(T, T') \le \sum_{k=2}^J \mathrm{diam}(T_k)$. In case (B), we note the Proposition 9 (ii) predicts that $\mathrm{gen}(T_{J+1}) = \mathrm{gen}(T)+1$ and that $T_{J+1}$ and $T_J$ are compatibly divisible. Hence, $\widetilde{T} = T_{J+1}$, and $\widetilde{T}$ is bisected only once. Therefore, the sons $T'$ of $\widetilde{T}$ fulfil $\mathrm{gen}(T') = \mathrm{gen}(\widetilde{T}) + 1 = \mathrm{gen}(T) + 2$ and $\mathrm{dist}(T', T_J) = 0$, which gives $\mathrm{dist}(T, T') \le \sum_{k=2}^J \mathrm{diam}(T_k)$. In the case (A) or (B), we can therefore estimate

$$\begin{aligned}\mathrm{dist}(T, T') &\le \sum_{k=2}^{J} \mathrm{diam}(T_k) \le C^{\mathrm{diam}} \sum_{k=2}^{J} 2^{-\mathrm{gen}(T_k)/2} \le C_{\mathrm{chain}} C^{\mathrm{diam}} 2^{-\mathrm{gen}(T)/2} \\ &= 2C_{\mathrm{chain}} C^{\mathrm{diam}} 2^{-(\mathrm{gen}(T)+2)/2} \le C_{\mathrm{dist}} C^{\mathrm{diam}} \sum_{k=\mathrm{gen}(T')}^{\mathrm{gen}(T)+2} 2^{-k/2}.\end{aligned}$$

It remains to consider the case (C). We use the triangle inequality to see

$$\mathrm{dist}(T, T') \le \mathrm{dist}(T_{J+1}, T') + \mathrm{diam}(T_{J+1}) + \sum_{k=2}^{J} \mathrm{diam}(T_k).$$

With $\mathrm{gen}(T_{J+1}) < \mathrm{gen}(T)$, Proposition 9 (ii) implies $\mathrm{gen}(T_{J+1}) = \mathrm{gen}(T_J) - i$ for $i \in \{1, 2\}$. We can use the induction hypothesis to estimate the first term on the right-hand side.

Furthermore, $\operatorname{diam}(T_{J+1}) \le C^{\mathrm{diam}} 2^{-(\operatorname{gen}(T)-i)/2}$ and $\operatorname{diam}(T_k) \le C^{\mathrm{diam}} 2^{-\operatorname{gen}(T)/2}$ show

$$\begin{aligned}
\operatorname{dist}(T,T') &\le C_{\mathrm{dist}} C^{\mathrm{diam}} \sum_{k=\operatorname{gen}(T')}^{\operatorname{gen}(T)-i+2} 2^{-k/2} + C^{\mathrm{diam}} 2^{-(\operatorname{gen}(T)-i)/2} + C_{\mathrm{chain}} C^{\mathrm{diam}} 2^{-\operatorname{gen}(T)/2} \\
&= C_{\mathrm{dist}} C^{\mathrm{diam}} \sum_{k=\operatorname{gen}(T')}^{\operatorname{gen}(T)-i+2} 2^{-k/2} + C^{\mathrm{diam}} 2^{-\operatorname{gen}(T)/2} 2^{i/2} + 2C_{\mathrm{chain}} C^{\mathrm{diam}} 2^{-(\operatorname{gen}(T)+2)/2}. \qquad (15)
\end{aligned}$$

In case of $i = 1$, Equation (15) yields

$$\begin{aligned}
\operatorname{dist}(T,T') &\le C_{\mathrm{dist}} C^{\mathrm{diam}} \sum_{k=\operatorname{gen}(T')}^{\operatorname{gen}(T)+1} 2^{-k/2} + C^{\mathrm{diam}} 2^{-(\operatorname{gen}(T)+2)/2} 2^{3/2} + 2C_{\mathrm{chain}} C^{\mathrm{diam}} 2^{-(\operatorname{gen}(T)+2)/2} \\
&\le C_{\mathrm{dist}} C^{\mathrm{diam}} \sum_{k=\operatorname{gen}(T')}^{\operatorname{gen}(T)+1} 2^{-k/2} + (2C_{\mathrm{chain}} + 2^{3/2}) C^{\mathrm{diam}} 2^{-(\operatorname{gen}(T)+2)/2} \\
&= C_{\mathrm{dist}} C^{\mathrm{diam}} \sum_{k=\operatorname{gen}(T')}^{\operatorname{gen}(T)+2} 2^{-k/2}.
\end{aligned}$$

In case of $i = 2$, Equation (15) yields

$$\begin{aligned}
\operatorname{dist}(T,T') &\le C_{\mathrm{dist}} C^{\mathrm{diam}} \sum_{k=\operatorname{gen}(T')}^{\operatorname{gen}(T)} 2^{-k/2} + C^{\mathrm{diam}} 2^{-(\operatorname{gen}(T)+1)/2} 2^{3/2} + 2C_{\mathrm{chain}} C^{\mathrm{diam}} 2^{-(\operatorname{gen}(T)+2)/2} \\
&\le C_{\mathrm{dist}} C^{\mathrm{diam}} \sum_{k=\operatorname{gen}(T')}^{\operatorname{gen}(T)} 2^{-k/2} + C_{\mathrm{dist}} C^{\mathrm{diam}} 2^{-(\operatorname{gen}(T)+1)/2} + C_{\mathrm{dist}} C^{\mathrm{diam}} 2^{-(\operatorname{gen}(T)+2)/2} \\
&= C_{\mathrm{dist}} C^{\mathrm{diam}} \sum_{k=\operatorname{gen}(T')}^{\operatorname{gen}(T)+2} 2^{-k/2}.
\end{aligned}$$

Altogether, this concludes the proof. □

*Proof of Theorem 5.* The proof is divided into four steps: First, we show the estimate (4) for Algorithm 1, i.e. `refineNVB`. Bootstrapping this result, we show the estimate (4) for `refineNVB3`, i.e. Algorithm 2 with the usual NVB-refinement patterns. In the third step, we show estimate (4) for Algorithm `refineNVBred`. Finally, we can boostrap the last result to remove any restrictions and show (4) for `refine`.

- **Proof for `refineNVB`:** We define $\mathcal{M} := \bigcup_{j=0}^{k-1} \mathcal{M}_j$ and note that the $\mathcal{M}_j$ are pairwise disjoint. We choose two sequences of positive, real numbers $(a_\ell)_{\ell=-2}^{\infty}$ und $(b_\ell)_{\ell=0}^{\infty}$ with the properties

$$\sum_{\ell=-2}^{\infty} a_\ell =: A < \infty, \qquad \sum_{\ell=0}^{\infty} 2^{-\ell/2} b_\ell =: B < \infty, \qquad \inf_{\ell \ge 1} b_\ell a_\ell =: c > 0$$

and $b_0 \geq 1$. A valid choice is given in Remark 19 below. Define a function $\lambda : \mathcal{T}_k \times \mathcal{M} \to \mathbb{R}_+$ as

$$\lambda(T, T_{\mathcal{M}}) = \begin{cases} a_{\mathrm{gen}(T_{\mathcal{M}})-\mathrm{gen}(T)}, & \text{if } \mathrm{dist}(T, T_{\mathcal{M}}) < C_\lambda B 2^{-\mathrm{gen}(T)/2} \text{ and } \mathrm{gen}(T) \leq \mathrm{gen}(T_{\mathcal{M}}) + 2, \\ 0 & \text{otherwise.} \end{cases}$$

Here, the constant $C_\lambda$ is defined as $C_\lambda := C^{\mathrm{diam}} \left(1 + \frac{C_{\mathrm{dist}}}{1-2^{-1/2}}\right)$. In the following, we show

$$\sum_{T \in \mathcal{T}_k \setminus \mathcal{T}_0} \lambda(T, T_{\mathcal{M}}) \leq C^{\$} \quad \text{for all } T_{\mathcal{M}} \in \mathcal{M} \text{ and} \tag{16}$$

$$\sum_{T_{\mathcal{M}} \in \mathcal{M}} \lambda(T, T_{\mathcal{M}}) \geq C_{\$} \quad \text{for all } T \in \mathcal{T}_k \setminus \mathcal{T}_0 \tag{17}$$

with some constants $C_{\$}$, $C^{\$} > 0$. The desired bound then follows from the last two estimates:

$$\begin{aligned} C_{\$} \left(\#\mathcal{T}_k - \#\mathcal{T}_0\right) \leq C_{\$} \#(\mathcal{T}_k \setminus \mathcal{T}_0) \leq \sum_{T \in \mathcal{T}_k \setminus \mathcal{T}_0} \sum_{T_{\mathcal{M}} \in \mathcal{M}} \lambda(T, T_{\mathcal{M}}) &= \sum_{T_{\mathcal{M}} \in \mathcal{M}} \sum_{T \in \mathcal{T}_k \setminus \mathcal{T}_0} \lambda(T, T_{\mathcal{M}}) \\ &\leq C^{\$} \#\mathcal{M} = C^{\$} \sum_{j=0}^{k-1} \#\mathcal{M}_j. \end{aligned}$$

We start with the proof of (16): Choose $T_{\mathcal{M}} \in \mathcal{M}$. For $g \in \mathbb{N}_0$ with $0 \leq g \leq \mathrm{gen}(T_{\mathcal{M}}) + 2$ we choose a subset of $\mathcal{T}_k$

$$\mathcal{T}_k(T_{\mathcal{M}}, g) := \left\{ T \in \mathcal{T}_k \mid \mathrm{dist}(T, T_{\mathcal{M}}) \leq C_\lambda B 2^{-g/2} \text{ and } \mathrm{gen}(T) = g \right\}.$$

For $0 \leq g \leq \mathrm{gen}(T_{\mathcal{M}}) + 2$ holds $\mathrm{diam}(T_{\mathcal{M}}) \leq C^{\mathrm{diam}} 2^{-\mathrm{gen}(T_{\mathcal{M}})/2} \leq 2C^{\mathrm{diam}} 2^{-g/2}$, and for $T \in \mathcal{T}_k(T_{\mathcal{M}}, g)$ holds $\mathrm{diam}(T) \leq C^{\mathrm{diam}} 2^{-g/2}$. Therefore, $\bigcup \mathcal{T}_k(T_{\mathcal{M}}, g)$ is contained in a ball with radius $r \simeq 2^{-g/2}$ around the center of mass of $T_{\mathcal{M}}$. For $T \in \mathcal{T}_k(T_{\mathcal{M}}, g)$ holds $|T| \geq C_{\mathrm{diam}}^2 2^{-g}$, hence the number of elements in $\mathcal{T}_k(T_{\mathcal{M}}, g)$ is bounded uniformly,

$$\#\mathcal{T}_k(T_{\mathcal{M}}, g) \leq C_4$$

with some constant $C_4 > 0$. We conclude

$$\sum_{T \in \mathcal{T}_k \setminus \mathcal{T}_0} \lambda(T, T_{\mathcal{M}}) = \sum_{g=0}^{\mathrm{gen}(T_{\mathcal{M}})+2} \sum_{T \in \mathcal{T}_k(T_{\mathcal{M}}, g)} a_{\mathrm{gen}(T_{\mathcal{M}})-g} \leq C_4 \sum_{\ell=-2}^{\infty} a_\ell = C_4 A.$$

It remains to show (17): We choose an arbitrary element $T_0 \in \mathcal{T}_k \setminus \mathcal{T}_0$. Then, there is an index $1 \leq k_0 \leq k$ such that $T_0 \in \mathcal{T}_{k_0}$ and $T_0 \notin \mathcal{T}_{k_0-1}$. Hence, there has to be an element $T_1 \in \mathcal{M}_{k_0-1} \subset \mathcal{M}$ such that $T_0$ is created by bisection of $T_1$. In an iterative manner, we can construct a sequence $(T_j)_{j=1}^{J} \subset \mathcal{M}$ such that $T_j$ is created by bisection of $T_{j+1}$.

According to Lemma 18 holds $\mathrm{gen}(T_{j+1}) \geq \mathrm{gen}(T_j) - 2$. Because $T_J \in \mathcal{T}_0$, we may choose the smallest index $s \in \mathbb{N}$ with $\mathrm{gen}(T_s) < \mathrm{gen}(T_0)$ and note $\mathrm{gen}(T_0) - 2 \leq \mathrm{gen}(T_s)$.

For $1 \leq j \leq s$ and $\ell > 0$ we define the set

$$\mathcal{T}_k(T_0, i, j) := \{T \in \{T_0, \ldots, T_{j-1}\} \mid \mathrm{gen}(T) = \mathrm{gen}(T_0) + i\}$$

with $m(i,j) := \#\mathcal{T}_k(T_0,i,j)$. We have

$$\operatorname{dist}(T_0,T_j) \le \sum_{i=1}^{j} \operatorname{dist}(T_{i-1},T_i) + \sum_{i=1}^{j-1} \operatorname{diam}(T_i).$$

According to Lemma 18,

$$\begin{aligned}
\operatorname{dist}(T_0,T_j) &\le \frac{C_{\text{dist}}C^{\text{diam}}}{1-2^{-1/2}} \sum_{i=1}^{j} 2^{-\operatorname{gen}(T_{i-1})/2} + C^{\text{diam}} \sum_{i=1}^{j-1} 2^{-\operatorname{gen}(T_i)/2} \\
&\le C^{\text{diam}} \left(1 + \frac{C_{\text{dist}}}{1-2^{-1/2}}\right) \sum_{i=0}^{j-1} 2^{-\operatorname{gen}(T_i)/2} \\
&= C_\lambda \sum_{i=0}^{j-1} 2^{-\operatorname{gen}(T_i)/2} \\
&= C_\lambda \sum_{i=0}^{\infty} m(i,j) 2^{-(\operatorname{gen}(T_0)+i)/2}.
\end{aligned} \tag{18}$$

We distinguish two cases:

(i) $m(i,s) \le b_i$ for all $i \in \mathbb{N}$. According to (18),

$$\operatorname{dist}(T_0,T_s) \le C_\lambda 2^{-\operatorname{gen}(T_0)/2} \sum_{i=0}^{\infty} b_i 2^{-i/2} = C_\lambda B 2^{-\operatorname{gen}(T_0)/2}.$$

Clearly,

$$\sum_{T_\mathcal{M} \in \mathcal{M}} \lambda(T_0, T_\mathcal{M}) \ge \lambda(T_0,T_s),$$

and the choice of $s$ provides $\operatorname{gen}(T_s) - \operatorname{gen}(T_0) \in \{-1,-2\}$. The definition of $\lambda$ shows

$$\lambda(T_0,T_s) = a_{\operatorname{gen}(T_s)-\operatorname{gen}(T_0)} \ge \min(a_{-2},a_{-1}) > 0.$$

(ii) There is $i \ge 0$ with $m(i,s) > b_i$. For all those $i$ exists a smallest $s(i)$ such that $m(i,s(i)) > b_i$. Choose $i^\star \in \underset{\substack{i\in\mathbb{N} \\ m(i,s(i))>b_i}}{\operatorname{argmin}}\, s(i)$. In particular, $1 \le s(i^\star) \le s$, and monotonicity of $m(i,j)$ in $j$ then shows

$$m(i,s(i^\star)-1) \le b_i \quad \text{for all } i \ge 0. \tag{19}$$

Furthermore,

$$m(i^\star, s(i^\star)) > b_{i^\star} \tag{20}$$

according to the choice of $i^\star$. Given $i \in \mathbb{N}$, equation (19) implies $\operatorname{dist}(T_0,T) \le C_\lambda B 2^{-\operatorname{gen}(T_0)/2}$ for $T \in \mathcal{T}_k(T_0,i,s(i^\star))$ as in (i). For $i^\star = 0$, at least $T_0 \in \mathcal{T}_k(T_0,0,s(0))$. However, equation (20) reveals $m(0,s(0)) > b_0 \ge 1$, whence there exists $T' \in \mathcal{T}_k(T_0,0,s(0)) \cap \mathcal{M}$. We conclude

$$\sum_{T_\mathcal{M}\in\mathcal{M}} \lambda(T_0,T_\mathcal{M}) \ge \lambda(T_0,T') = a_0 > 0.$$

Now suppose $i^\star > 0$. Clearly, $T_0 \notin \mathcal{T}_k(T_0, i^\star, s(i^\star))$ whence $\mathcal{T}_k(T_0, i^\star, s(i^\star)) \subset \mathcal{M}$. Per definition, $\lambda(T_0, T') = a_{i^\star}$ for all $T' \in \mathcal{T}_k(T_0, i^\star, s(i^\star))$. This reveals

$$\sum_{T_\mathcal{M}\in\mathcal{M}} \lambda(T_0, T_\mathcal{M}) \geq \sum_{T_\mathcal{M}\in\mathcal{T}_k(T_0,i^\star,s(i^\star))} \lambda(T_0, T_\mathcal{M}) = m(i^\star, s(i^\star))a_{i^\star}$$
$$> b_{i^\star} a_{i^\star} \geq \inf_{i\geq 1} b_i a_i = c > 0.$$

Setting $C_\$ := \min\{a_{-2}, a_{-1}, a_0, c\}$ shows (17) and thus concludes the proof for algorithm `refineNVB`.

- **Proof for `refineNVB3`:** As stated in Remark 3, a call of `refineNVB3` can be carried out by two calls of `refineNVB`, hence

$$\#\mathcal{T}_\ell - \#\mathcal{T}_0 \lesssim \sum_{j=0}^{\ell-1} \#\mathcal{M}_j + \sum_{j=0}^{\ell-1} \#\mathcal{M}_{j+1/2}.$$

Since $\#\mathcal{M}_{j+1/2} \leq 2\#\mathcal{M}_j$, we conclude

$$\#\mathcal{T}_\ell - \#\mathcal{T}_0 \lesssim \sum_{j=0}^{\ell-1} \#\mathcal{M}_j.$$

- **Proof for `refineNVBred`:** We consider the meshes $\widetilde{\mathcal{T}}_\ell$ that are constructed in Lemma 14. Using the notation from Lemma 14 and Proposition 13, we conclude

$$\#\mathcal{T}_\ell - \#\mathcal{T}_0 = \#\widetilde{\mathcal{T}}_\ell - \#\widetilde{\mathcal{T}}_0 \lesssim \sum_{j=0}^{\ell-1} \#\widetilde{\mathcal{M}}_j \lesssim \sum_{j=0}^{\ell-1} \#\mathcal{M}_j,$$

since $\#\widetilde{\mathcal{M}}_j \leq 2\#\mathcal{M}_j$.

- **Proof for `refine`:** A step of `refine` can be carried out by two steps of `refineNVBred`, and we may argue as in the last steps.

□

**Remark 19.** *A valid choice for $(a_\ell)_{\ell=-2}^\infty$ and $(b_\ell)_{\ell=0}^\infty$ is $b_\ell = 2^{\ell/3}$ and $a_{-2} = 1$, $a_\ell = (\ell+2)^{-2}$ for $\ell \geq -1$. In this case, $A = 1 + \pi^2/6$ and $B = (1 - 2^{-1/6})^{-1}$.*

## 5. Proof of Theorem 6 ($H^1$-stability of $L_2$-projection)

To prove Theorem 6, we mainly follow the ideas of [7]: a deliberate choice of nodal values for a mesh $\mathcal{T}_\ell$ (Proposition 20) yields local estimate for the eigenvalues of the mass matrices (Lemma 22), from which $H^1$-stability of $\Pi_\ell$ follows readily. These ideas are already used in [7]. Our choice of the nodal values is also inspired by the latter work, but the analysis clearly needs certain refined arguments to adapt it to the present situation.

Throughout this section, we use the following notation: Each element $T \in \mathcal{T}_\ell$ is the convex hull of its three nodes $T = \text{conv}\{z_{[T,1]}, z_{[T,2]}, z_{[T,3]}\}$ with certain nodes $z_{[T,i]} \in \mathcal{N}_\ell$. With these local coordinates, we define the symmetric element mass matrix

$$\mathbf{M}_T \in \mathbb{R}^{3\times 3}, \quad (\mathbf{M}_T)_{ij} = \int_T \varphi_{[T,i]}(\mathbf{x})\varphi_{[T,j]}(\mathbf{x})\, d\mathbf{x} \quad \text{for all } T \in \mathcal{T}_\ell \text{ and } i,j = 1,2,3.$$

Here $\varphi_k \in \mathcal{S}^1(\mathcal{T}_\ell)$ denotes the hat function associated with the node $z_k \in \mathcal{N}_\ell$. Moreover, for each node $z_k \in \mathcal{N}_k$, let $d_k > 0$ denote some quantity which is specified below and satisfies

$d_{[T,i]} \simeq |T|^{1/2} \simeq \operatorname{diam}(T) =: h_\ell(T)$ for all $T \in \mathcal{T}_\ell$ and $i = 1, 2, 3$. We define the diagonal element scaling matrix

$$(\Lambda_T) \in \mathbb{R}^{3\times 3}, \quad (\Lambda_T)_{ij} = \frac{h_\ell(T)}{d_{[T,i]}}\,\delta_{ij} \quad \text{for all } T \in \mathcal{T}_\ell \text{ and } i, j = 1, 2, 3,$$

where $\delta_{ij} \in \{0, 1\}$ denotes Kronecker's delta. Besides the precise definition of the quantities $d_k > 0$, the essential step in the proof of Theorem 6 is stated in the following proposition.

**Proposition 20.** *Let $\mathcal{T}_\ell$ be a $\gamma$-shape regular mesh on $\Omega$. Suppose that for every $T \in \mathcal{T}_\ell$ holds*

$$\frac{d_{[T,i]}}{d_{[T,j]}} \le C_5 \le \pi \quad \textit{for all } i, j = 1, 2, 3 \tag{21}$$

*as well as*

$$C_6^{-1} \le \frac{d_{[T,i]}}{h_\ell(T)} \le C_6 \quad \textit{for all } i = 1, 2, 3 \tag{22}$$

*with constants $C_5, C_6 > 0$. Then, the $L_2$-projection $\Pi_\ell$ onto $\mathcal{S}^1(\mathcal{T}_\ell)$ is $H^1$-stable* (9)*, and the constant $C_2 > 0$ depends only on $C_5, C_6$, and the $\gamma$-shape regularity of $\mathcal{T}_\ell$.*

The proof requires the following elementary estimate.

**Lemma 21.** *Suppose that $a, b, c > 0$ and $M \ge 1$ are real numbers with $M^{-1} \le a, b, c \le M$ and $c = ab$. Then,*

$$a + b + c + a^{-1} + b^{-1} + c^{-1} \le 2(1 + M + M^{-1}),$$

*and the upper bound in this estimate is sharp.* □

The next lemma provides a criterion for $H^1$-stability of the $L_2$-projection in terms of local (generalized) eigenvalue estimates.

**Lemma 22.** *Let $\mathcal{T}_\ell$ be a $\gamma$-shape-regular mesh on $\Omega$ such that for all $T \in \mathcal{T}_\ell$ holds*

$$C_7^{-1}\mathbf{x}^\top \Lambda_T^2 \mathbf{M}_T \Lambda_T^2 \mathbf{x} \le \mathbf{x}^\top \mathbf{M}_T \mathbf{x} \le C_8 \mathbf{x}^\top \Lambda_T^2 \mathbf{M}_T \mathbf{x} \quad \textit{for all } \mathbf{x} \in \mathbb{R}^3 \tag{23}$$

*with constants $C_7, C_8 > 0$ which do not depend on $T \in \mathcal{T}_\ell$. Then, the $L_2$-projection $\Pi_\ell$ onto $\mathcal{S}^1(\mathcal{T}_\ell)$ is $H^1$-stable* (9)*, and the constant $C_2 > 0$ depends only on $\gamma$-shape regularity and on $C_7, C_8 > 0$.*

*Proof.* Let $J_\ell : H^1(\Omega) \to \mathcal{S}^1(\mathcal{T}_\ell)$ denote a Clément-type quasi interpolation operator, i.e. it holds that

$$\|\nabla J_\ell u\|_{L_2(T)} \lesssim \|\nabla u\|_{L_2(\omega_T)} \quad \text{as well as} \quad \|u - J_\ell u\|_{L_2(T)} \lesssim \|h_\ell \nabla u\|_{L_2(\omega_T)}, \tag{24}$$

for all $u \in H^1(\Omega)$ and some hidden constant $C > 0$ which depends only on $\gamma$-shape regularity. Here, $h_\ell \in L^\infty(\Omega)$ denotes the $\mathcal{T}_\ell$-piecewise constant mesh-size function $h_\ell|_T = h_\ell(T)$ and $\omega_T := \bigcup \{T' \in \mathcal{T}_\ell : T' \cap T \neq \emptyset\}$ is the patch of elements around $T$. A standard local inverse estimate and stability of $J_\ell$ show

$$\|\nabla \Pi_\ell u\|_{L_2(\Omega)} \le \|\nabla(\Pi_\ell u - J_\ell u)\|_{L_2(\Omega)} + \|\nabla J_\ell u\|_{L_2(\Omega)} \lesssim \|h_\ell^{-1}(\Pi_\ell u - J_\ell u)\|_{L_2(\Omega)} + \|\nabla u\|_{L_2(\Omega)},$$

where the hidden constant depends only on $\gamma$-shape regularity. With $q_\ell := \Pi_\ell u - J_\ell u \in \mathcal{S}^1(\mathcal{T}_\ell)$, it thus only remains to bound $\|h_\ell^{-1} q_\ell\|_{L_2(\Omega)}$ by $\|\nabla u\|_{L_2(\Omega)}$. To that end, let $\mathcal{N}_\ell =$ $\{z_1, \ldots, z_n\}$ and define $p_\ell = \sum_{j=1}^n q_j d_j^{-2} \varphi_j \in \mathcal{S}^1(\mathcal{T}_\ell)$. With $\mathbf{x}_T = (q_{[T,1]}, q_{[T,2]}, q_{[T,3]})^\top$ and the upper estimate in (23), we see

$$\|h_\ell^{-1} q_\ell\|_{L_2(\Omega)}^2 = \sum_{T\in\mathcal{T}_\ell} h_\ell(T)^{-2} \mathbf{x}_T^\top \mathbf{M}_T \mathbf{x}_T \le C_8 \sum_{T\in\mathcal{T}_\ell} h_\ell(T)^{-2} \mathbf{x}_T^\top \Lambda_T^2 \mathbf{M}_T \mathbf{x}_T$$
$$= C_8 \sum_{T\in\mathcal{T}_\ell} \sum_{i=1}^3 \frac{q_{[T,i]}}{d_{[T,i]}^2} \int_T \varphi_{[T,i]} q_\ell = C_8 \int_\Omega p_\ell q_\ell.$$

Next, we use the orthogonality relation of $\Pi_\ell$ and the approximation property of $J_\ell$ to see

$$\int_\Omega p_\ell q_\ell = \int_\Omega p_\ell(\Pi_\ell u - J_\ell u) = \int_\Omega p_\ell(u - J_\ell u) \lesssim \|h_\ell p_\ell\|_{L_2(\Omega)} \|\nabla u\|_{L_2(\Omega)}.$$

The combination of the last two estimates gives

$$\|h_\ell^{-1} q_\ell\|_{L_2(\Omega)}^2 \lesssim C_8\, \|h_\ell p_\ell\|_{L_2(\Omega)} \|\nabla u\|_{L_2(\Omega)}.$$

With $\mathbf{y}_T = (q_{[T,1]} d_{[T,1]}^{-2}, q_{[T,2]} d_{[T,2]}^{-2}, q_{[T,3]} d_{[T,3]}^{-2})^\top = h_\ell(T)^{-2}\, \Lambda_T \mathbf{x}_T$ and the lower estimate in (23), we obtain

$$\|h_\ell p_\ell\|_{L_2(\Omega)}^2 = \sum_{T\in\mathcal{T}_\ell} h_\ell(T)^2\, \mathbf{y}_T^\top \mathbf{M}_T \mathbf{y}_T = \sum_{T\in\mathcal{T}_\ell} h_\ell(T)^{-2} \mathbf{x}_T^\top \Lambda_T^2 \mathbf{M}_T \Lambda_T^2 \mathbf{x}_T$$
$$\le C_7 \sum_{T\in\mathcal{T}_\ell} h_\ell(T)^{-2} \mathbf{x}_T^\top \mathbf{M}_T \mathbf{x}_T = C_7\, \|h_\ell^{-1} q_\ell\|_{L_2(\Omega)}^2.$$

The combination of the last two estimates reveals

$$\|h_\ell^{-1} q_\ell\|_{L_2(\Omega)}^2 \lesssim C_7 C_8\, \|h_\ell^{-1} q_\ell\|_{L_2(\Omega)} \|\nabla u\|_{L_2(\Omega)}.$$

Division by $\|h_\ell^{-1} q_\ell\|_{L_2(\Omega)}$ thus concludes the proof. □

*Proof of Proposition 20.* We will verify the eigenvalue criterion (23) of Lemma 22: Let $\mathbf{x} \in \mathbb{R}^3$. For each element $T \in \mathcal{T}_\ell$, a standard scaling argument proves equivalence

$$|T|\, \mathbf{x}^\top \mathbf{x} \lesssim \mathbf{x}^\top \mathbf{M}_T \mathbf{x} \lesssim |T|\, \mathbf{x}^\top \mathbf{x}, \tag{25}$$

cf. e.g. [17, Proposition 6.3.1], where the hidden constants depend only on $\gamma$-shape regularity. With (25) and the lower bound of Assumption (22) we obtain

$$\mathbf{x}^\top \Lambda_T^2 \mathbf{M}_T \Lambda_T^2 \mathbf{x} \lesssim |T| (\Lambda_T^2 \mathbf{x})^\top (\Lambda_T^2 \mathbf{x}) \le C_6^4\, |T|\, \mathbf{x}^\top \mathbf{x} \lesssim C_6^4\, \mathbf{x}^\top \mathbf{M}_T \mathbf{x}.$$

This proves the lower estimate in (23), and $C_7$ depends only on $C_6$ and $\gamma$-shape regularity.

The verification of the upper estimate in (23) is more involved: We fix a reference element $\widehat{T}$ such that the corresponding element mass matrix satisfies $\widehat{\mathbf{M}}_{jk} = 1 + \delta_{jk}$ with Kronecker's delta. The transformation theorem yields $\mathbf{M}_T = C_{\mathrm{ref}} |T|\, \widehat{\mathbf{M}}$. We define the matrix

$$\mathbf{A}_T := \Lambda_T^2 \mathbf{M}_T + \mathbf{M}_T \Lambda_T^2 = C_{\mathrm{ref}} |T| \left(\Lambda_T^2 \widehat{\mathbf{M}} + \widehat{\mathbf{M}} \Lambda_T^2\right) =: C_{\mathrm{ref}} |T|\, \widehat{\mathbf{A}}_T.$$

Symmetry of $\Lambda_T$ and $\mathbf{M}_T$ gives $\mathbf{x}_T^\top \mathbf{M}_T \Lambda_T^2 \mathbf{x} = (\Lambda_T^2 \mathbf{x})^\top \mathbf{M}_T \mathbf{x} = \mathbf{x}_T^\top \Lambda_T^2 \mathbf{M}_T \mathbf{x}$, whence

$$\mathbf{x}^\top \mathbf{A}_T \mathbf{x} = 2\, \mathbf{x}^\top \Lambda_T^2 \mathbf{M}_T \mathbf{x}. \tag{26}$$

Moreover, it holds that

$$(\widehat{\mathbf{B}}_T)_{jk} := \left(\Lambda_T^{-1} \widehat{\mathbf{A}}_T \Lambda_T^{-1}\right)_{jk} = \left(\Lambda_T \widehat{\mathbf{M}} \Lambda_T^{-1} + \Lambda_T^{-1} \widehat{\mathbf{M}} \Lambda_T\right)_{jk} = \left(\frac{d_{[T,j]}}{d_{[T,k]}} + \frac{d_{[T,k]}}{d_{[T,j]}}\right) \widehat{\mathbf{M}}_{jk}.$$

In particular, the symmetric matrix $\widehat{\mathbf{B}}_T \in \mathbb{R}^{3\times 3}$ satisfies $(\widehat{\mathbf{B}}_T)_{ii} = 4$ and $(\widehat{\mathbf{B}}_T)_{ij} = b_i/b_j + b_j/b_i$ with certain $b_i = d_{[T,i]} > 0$ for all $i, j = 1, 2, 3$. We apply [4, Proposition 6.1] and obtain that the minimal eigenvalue $\lambda_{\min}$ of $\widehat{\mathbf{B}}_T$ satisfies

$$\lambda_{\min} = 5 - \Big(\sum_{j,k=1}^{3} \frac{d_{[T,j]}^2}{d_{[T,k]}^2}\Big)^{1/2}$$

Now, we employ Lemma 21 with $a = d_{[T,1]}/d_{[T,2]}$, $b = d_{[T,2]}/d_{[T,3]}$, and $c = d_{[T,1]}/d_{[T,3]} = ab$ as well as $M = C_5$. This gives

$$\sum_{j,k=1}^{3} \frac{d_{[T,j]}^2}{d_{[T,k]}^2} \le 3 + 2(1 + \pi^2 + \pi^{-2}) \le 24.95 < 25. \tag{27}$$

Consequently, $\lambda_{\min} > 0$ and thus $\widehat{\mathbf{B}}_T$ is positive definite with

$$\mathbf{y}^\top \widehat{\mathbf{B}}_T \mathbf{y} \ge \lambda_{\min}\, \mathbf{y}^\top \mathbf{y} \quad \text{for all } \mathbf{y} \in \mathbb{R}^3$$

We stress that $\lambda_{\min}$ depends only on $C_5 \le \pi$. Next, we consider $\mathbf{y} = \Lambda_T \mathbf{x}$ and note that, due to (22),

$$\mathbf{y}^\top \mathbf{y} \simeq \mathbf{x}^\top \mathbf{x},$$

where the hidden constants depend only on $C_6$. We combine the last two estimates to see

$$\mathbf{x}^\top \mathbf{M}_T \mathbf{x} \simeq |T|\, \mathbf{x}^\top \mathbf{x} \simeq |T|\, \mathbf{y}^\top \mathbf{y} \lesssim |T|\, \mathbf{y}^\top \widehat{\mathbf{B}}_T \mathbf{y} = |T|\, \mathbf{x}^\top \widehat{\mathbf{A}}_T \mathbf{x} \simeq \mathbf{x}^\top \mathbf{A}_T \mathbf{x} = 2\, \mathbf{x}^\top \Lambda_T^2 \mathbf{M}_T \mathbf{x},$$

where the final equality has been noted in (26). This is the upper bound in (23), and $C_8$ depends only on $C_5, C_6$ and on $\gamma$-shape regularity of $\mathcal{T}_\ell$. □

**Remark 23.** *As follows from* (27)*, the assumption $C_5 \le \pi$ from* (21) *can slightly be relaxed to $1 + C_5^2 + C_5^{-2} < 11$.*

*Proof of Theorem 6.* It only remains to define $d_j > 0$ for all $z_j \in \mathcal{N}_\ell$ such that the assumptions (21)–(22) of Proposition 20 are satisfied. To that end, we define $\delta(z_j, z_k) \in \mathbb{N}_0$ for all nodes $z_j, z_k \in \mathcal{N}_\ell$ as follows:

- $\delta(z_j, z_k) := 0$ if $z_j = z_k$,
- $\delta(z_j, z_k) := 1$ if there exist some $T \in \mathcal{T}_\ell$ with $z_j, z_k \in T$,
- $\delta(z_j, z_k) := N$ if $N \ge 2$ is the minimal number of elements $T_1, \dots, T_N \in \mathcal{T}_\ell$ with $z_j \in T_1$, $T_j \cap T_{j+1} \in \mathcal{E}_\ell$ for all $j = 1, \dots, N-1$, and $z_k \in T_N$.

Finally, we now define

$$d_j := \min_{T \in \mathcal{T}_\ell} 2^{(2\delta(z_j,T) - \mathrm{gen}(T))/2}, \quad \text{where} \quad \delta(z_j, T) := \min_{z_k \in \mathcal{N}_\ell \cap T} \delta(z_j, z_k). \tag{28}$$

● We first prove (21): Let $T \in \mathcal{T}_\ell$ and $z_j, z_k \in \mathcal{N}_\ell \cap T$. Suppose that $T' \in \mathcal{T}_\ell$ satisfies $d_k = 2^{(2\delta(z_k,T') - \mathrm{gen}(T'))/2}$. By definition, it holds that $d_j \le 2^{(2\delta(z_j,T') - \mathrm{gen}(T'))/2}$ and $|\delta(z_k, T') - \delta(z_j, T')| \le 1$. This implies

$$\frac{d_j}{d_k} \le \frac{2^{(2\delta(z_j,T') - \mathrm{gen}(T'))/2}}{2^{(2\delta(z_k,T') - \mathrm{gen}(T'))/2}} \le 2^{(\delta(z_j,T') - \delta(z_k,T'))} \le 2 \le \pi.$$

● Second, we prove the upper bound $d_j \le C_{\rm diam}^{-1} h_\ell(T)$ for all $T \in \mathcal{T}_\ell$ and $z_j \in \mathcal{N}_\ell \cap T$: By definition, it holds that $\delta(z_j, T) = 0$. With Proposition 9 (i), we thus see

$$d_j \le 2^{-\mathrm{gen}(T/2)} \le C_{\rm diam}^{-1} |T|^{1/2} \le C_{\rm diam}^{-1}\, h_\ell(T).$$

● Third, we prove the lower bound in (22) for all $T \in \mathcal{T}_\ell$ and $z_j \in \mathcal{N}_\ell \cap T$: As above, we choose $T' \in \mathcal{T}_\ell$ with $d_j = 2^{(2\delta(z_j,T') - \mathrm{gen}(T'))/2}$. Note that $\mathrm{gen}(T') < \mathrm{gen}(T)$ would imply

$$2^{(2\delta(z_j,T)-\mathrm{gen}(T))/2} = 2^{-\mathrm{gen}(T)/2} < 2^{-\mathrm{gen}(T')/2} \le 2^{(2\delta(z_j,T')-\mathrm{gen}(T'))/2} = d_j$$

and hence a contradiction to the choice of $T'$. We thus see $\mathrm{gen}(T) \le \mathrm{gen}(T')$. This gives

$$\begin{aligned} d_j^2 = 2^{2\delta(z_j,T')-\mathrm{gen}(T')} &= 2^{-\mathrm{gen}(T)}\, 2^{2\delta(z_j,T')-|\mathrm{gen}(T')-\mathrm{gen}(T)|} \\ &\ge \big(C^{\rm diam}\big)^{-2}\, h_\ell(T)^2\, 2^{2\delta(z_j,T')-|\mathrm{gen}(T')-\mathrm{gen}(T)|} \end{aligned}$$

according to Proposition 9 (i). For $\delta(z_j, T') = N$, there is a node $z_k \in \mathcal{N}_\ell \cap T'$ and elements $T_1, \dots, T_N$ such that $z_j \in T_1$, $T_j \cap T_{j+1} \in \mathcal{E}_\ell$ for all $j = 1, \dots, N-1$, and $z_k \in T_N$. According to $\gamma$-shape regularity, the number of elements having $z_j$ or $z_k$ as one of their nodes is uniformly bounded by some constant $C_{\rm patch} > 0$ which depends only on $\gamma$. Hence, we can find a sequence of elements $\widetilde{T}_1, \dots, \widetilde{T}_n$ with $n \le N + 2C_{\rm patch}$, $\widetilde{T}_1 = T$, $\widetilde{T}_j \cap \widetilde{T}_{j+1} \in \mathcal{E}_\ell$ for all $j = 1, \dots, n-1$, and $\widetilde{T}_n = T'$. Furthermore, Proposition 9 (ii) states $|\mathrm{gen}(\widetilde{T}_j) - \mathrm{gen}(\widetilde{T}_{j+1})| \le 2$. This yields

$$|\mathrm{gen}(T) - \mathrm{gen}(T')| \le \sum_{j=1}^{n-1} |\mathrm{gen}(\widetilde{T}_j) - \mathrm{gen}(\widetilde{T}_{j+1})| \le 2(n-1).$$

Finally, we combine this with $\delta(z_j, T') = N$ and $n \le N + 2C_{\rm patch}$ to see

$$2^{2\delta(z_j,T')-|\mathrm{gen}(T')-\mathrm{gen}(T)|} \ge 2^{2N-2(n-1)} \ge 2^{2-4C_{\rm patch}}.$$

From this we infer

$$d_j^2 \ge \big(C^{\rm diam}\big)^{-2}\, 2^{2-4C_{\rm patch}} h_\ell(T)^2.$$

● Altogether, the definition $C_6 := \max\Big\{C_{\rm diam}^{-1}, \big(C^{\rm diam}\big)^{-1} 2^{1-2C_{\rm patch}}\Big\}$ shows (22). We stress that the constant $C_6^{-1}$ depends only on $\mathcal{T}_0$ and the uniform $\gamma$-shape regularity of meshes $\mathcal{T}_\ell$ generated by MNVB. Now that we have checked the assumptions (21)–(22), Proposition 20 applies and proves $H^1$-stability (9) of the $L_2$-projection $\Pi_\ell$ onto $\mathcal{S}^1(\mathcal{T}_\ell)$. The constant $C_2 > 0$ depends only on the initial mesh $\mathcal{T}_0$. □

Institute for Analysis and Scientific Computing,Vienna University of Technology, Wiedner Hauptstrasse 8-10, A-1040 Wien, Austria

*E-mail address*: {David.Pavlicek, Dirk.Praetorius}@tuwien.ac.at

*E-mail address*: Michael.Karkulik@tuwien.ac.at